\documentclass{amsart}
\usepackage{amssymb}
\usepackage[v2,cmtip]{xy}

\chardef\bs=`\\ 
\newcommand{\ssdot}{\bullet}
\newcommand{\subdot}{_\ssdot}

\mathchardef\varSigma="7106
\mathchardef\varOmega="710A

\newcommand{\CAT}{\aC} 

\newcommand{\initial}{\emptyset} 

\newcommand{\Com}[2]{\aC_{#1}/#2} 
\newcommand{\hCom}[1]{\h[\Com{#1}{#1}]} 
\newcommand{\Mod}[1]{\aM_{#1}} 
\DeclareMathOperator\h{Ho} 

\newcommand{\TAQ}{D}
\newcommand{\cotan}[2]{\mathbf{L\varOmega}_{#1}#2}
\newcommand{\Bcotan}[2]{\mathbf{LAb}_{#1}^{B}#2} 

\newcommand{\newh}[2]{#1^{\mathrm{D}}_{#2}} 

\newcommand{\cone}[1]{\mathrm{C}}
\newcommand{\susp}[1]{\mathrm{E}}

\newcommand{\loops}[1]{\Omega_{#1}}

\newcommand{\spc}[1]{\mathrm{Sp}(#1)} 
\newcommand{\csp}[1]{\spectrumfont{#1}} 

\newcommand{\cspsrel}[1]{\susp{B}^{\infty}(B\sma_{R}#1)}

\newcommand{\csps}{\susp{B}^{\infty}} 

\newcommand{\sps}[1]{\susp{}^{\infty}{#1}}

\newcommand{\tpi}{\tilde\pi} 
\newcommand{\tsigma}{\tilde\sigma} 

\newcommand{\dsusp}[1]{\susp{#1}^{\dL}} 
\newcommand{\dloops}[1]{\loops{#1}^{\dR}} 
\newcommand{\dL}{\mathbf{L}}
\newcommand{\dR}{\mathbf{R}}

\newcommand{\rc}[2]{#1^{#2}} 
\newcommand{\rh}[2]{#1_{#2}} 
\newcommand{\sic}{\sigma} 
\newcommand{\sih}{\sigma} 

\newcommand{\NUCA}[1]{\aN_{#1}} 
\newcommand{\oNCom}{\oC\widetilde{om}} 
\newcommand{\bNCom}{\bN} 
\newcommand{\derbNCom}{\bN^{\dL}} 
\newcommand{\derI}{I^{\dR}} 
\newcommand{\derindec}[1]{Q^{\dL}} 
\newcommand{\indec}[1]{Q} 
\newcommand{\zero}[1]{Z} 

\newcommand{\stab}{\bS} 
\newcommand{\nstab}{S} 

\newcommand{\fancyzero}[1]{\mathrm{Z}#1} 

\newcommand{\owspcn}{W(\NUCA{B})}  
\newcommand{\uw}{W} 

\newcommand{\spapp}[1]{\csp{A}_{\csp{#1}}} 
\newcommand{\zspapp}[1]{A_{0}{\csp{#1}}}
\newcommand{\cofapp}[1]{M_{\csp{#1}}}
\newcommand{\bigisonat}{\phi}
\newcommand{\bigiso}[1]{\bigisonat_{\csp{#1}}}

\newcommand{\dcsps}{\susp{B}^{\dL\infty}} 
\newcommand{\cutoff}[2]{T_{#1}\csp{#2}}
\newcommand{\spfree}[1]{F_{#1}}
\newcommand{\restrict}[2]{#1^{\susp{}}_{#2}}
\newcommand{\extend}[2]{#1^{c}_{#2}}
\newcommand{\extres}[2]{#1^{c\susp{}}_{#2}}
\newcommand{\resext}[2]{#1^{\susp{}c}_{#2}}

\newcommand{\newss}{\Sigma^{\infty}_{S+}}
\newcommand{\newssp}{\Sigma^{\infty}_{S}}
\newcommand{\lsusp}{\mathrm{B}}
\newcommand{\dcspsS}{\susp{S}^{\dL\infty}} 
\newcommand{\scspsS}{\susp{S}^{\infty}} 
\newcommand{\lfree}{\mathrm{L}} 
\newcommand{\btensor}{\mathbin{\widehat \otimes}} 

\newcommand{\bCom}{\bP} 
\newcommand{\ISP}{I_{\mathrm{Sp}}}
\newcommand{\JSP}{J_{\mathrm{Sp}}}
\newcommand{\srs}[2]{\csp{S}{}_{#1;#2}}
\newcommand{\tar}[2]{\csp{T}{}_{#1;#2}}
\newcommand{\jmap}[2]{\lambda_{#1;#2}}
\newcommand{\ComB}{\aC_{B}} 
\newcommand{\Scotan}[1]{\mathbf{LAb}^{S}_{S}#1} 

\newcommand{\oCom}{\oC{om}} 
\newcommand{\oAss}{\oA{ss}} 
\newcommand{\oGen}{\oG} 
\newcommand{\bGen}{\bG} 
\newcommand{\aGen}{\aG_{R/B}} 
\newcommand{\aGenB}{\aG_{B/B}} 
\newcommand{\aGenN}{\aN_{\oNUn}} 
\newcommand{\op}{\mathrm{op}}
\newcommand{\oUn}{\oU B}
\newcommand{\oUnX}{\oU X}
\newcommand{\oNUn}{\tilde\oU B}
\newcommand{\bUn}{\bCom_{\oUn}}
\newcommand{\bNUn}{\bNCom_{\oNUn}}

\newcommand{\spectrumfont}[1]{\mathpalette\niceunderline{#1}#1}
\newcommand{\niceunderline}[2]{%
{\setbox0\hbox{$#1\mathrm{#2}$}\setbox2\hbox{$#1#2$}%
\ifdim\wd0<\wd2%
\hbox to 0pt{$\underline{\hbox to \wd0{\hss}}$\hss}%
\else%
\hbox to 0pt{$\underline{\hbox to \wd2{\hss}}$\hss}%
\fi}}

\newcommand{\bF}{{\mathbb{F}}}
\newcommand{\bG}{{\mathbb{G}}}

\newcommand{\bL}{{\mathbb{L}}}

\newcommand{\bN}{{\mathbb{N}}}

\newcommand{\bP}{{\mathbb{P}}}

\newcommand{\bS}{{\mathbb{S}}}
\newcommand{\bT}{{\mathbb{T}}}

\newcommand{\bZ}{{\mathbb{Z}}}

\let\catsymbfont\mathfrak

\newcommand{\aC}{{\catsymbfont{C}}}

\newcommand{\aG}{{\catsymbfont{G}}}

\newcommand{\aM}{{\catsymbfont{M}}}
\newcommand{\aN}{{\catsymbfont{N}}}

\let\opsymbfont\mathcal 

\newcommand{\oA}{{\opsymbfont{A}}}

\newcommand{\oC}{{\opsymbfont{C}}}

\newcommand{\oF}{{\opsymbfont{F}}}
\newcommand{\oG}{{\opsymbfont{G}}}

\newcommand{\oL}{{\opsymbfont{L}}}

\newcommand{\oU}{{\opsymbfont{U}}}

\newcommand{\iso}{\cong}     
\newcommand{\thp}{\ltimes}   
\newcommand{\sma}{\wedge}    

\newcommand{\einf}{\text{$E_\infty$} } 

\renewcommand{\to}{\mathchoice{\longrightarrow}{\rightarrow}{\rightarrow}{\rightarrow}}
\newcommand{\from}{\mathchoice{\longleftarrow}{\leftarrow}{\leftarrow}{\leftarrow}}

\newcommand{\tofrom}{\leftrightarrows}

\newcommand{\overto}[1]{\xrightarrow{#1}}
\newcommand{\overfrom}[1]{\xleftarrow{#1}}

\def\quickop#1{\expandafter\DeclareMathOperator\csname #1\endcsname{#1}}
\quickop{Hom} \quickop{End} \quickop{Aut} \quickop{Mic}
\quickop{Ext} \quickop{Tor} \quickop{Id} \quickop{id} \quickop{Coker}
\quickop{Ker} \quickop{Lim} \quickop{Colim} \quickop{Holim}
\quickop{Hocolim} \quickop{Cotor} \quickop{Der}\quickop{Tel}

\newtheorem{intthm}{Theorem}
\newtheorem*{intcor}{Corollary}
\newtheorem{thm}{Theorem}[section]
\newtheorem{cor}[thm]{Corollary}
\newtheorem{lem}[thm]{Lemma}
\newtheorem{prop}[thm]{Proposition}

\theoremstyle{definition}
\newtheorem*{intdefn}{Definition}
\newtheorem{defn}[thm]{Definition}

\theoremstyle{remark}
\newtheorem{rem}[thm]{Remark}

\numberwithin{equation}{section}

\makeatletter
\let\c@equation\c@thm
\makeatother

\renewcommand{\theenumi}{\roman{enumi}}\renewcommand{\labelenumi}{(\theenumi)}

\begin{document}

\title{Homology and Cohomology of \einf Ring Spectra}
\author{Maria Basterra}
\address{Department of Mathematics, University of New Hampshire, Durham, NH}
\email{basterra@math.unh.edu}
\author{Michael A. Mandell}
\address{Department of Mathematics, University of Chicago, Chicago, IL}
\email{mandell@math.uchicago.edu}
\thanks{Corresponding author. Fax: (603) 862-4096; E-mail address: basterra@math.unh.edu}
\thanks{The second author was supported in part by NSF grant DMS-0203980}

\date{June 30, 2004}

\begin{abstract}
Every homology or cohomology theory on a category of \einf
ring spectra is Topological Andr\'e--Quillen homology or cohomology with
appropriate coefficients.  Analogous results hold for the category of
$A_\infty$ ring spectra and for categories of algebras over many other
operads. 
\end{abstract}

\maketitle

\section*{Introduction}

Homology and cohomology theories in various contexts provide some of
the most effective tools in mathematics because of their
computability, their usefulness for classification problems,
and their close relationship to extensions and obstructions.  In the
context of homotopy theory, homology and cohomology typically refer to
theories that satisfy the Eilenberg--Steenrod (and Milnor) axioms.
Although originally phrased for topological spaces, these axioms make
sense in the more general context of homotopy theory of closed model
categories. 

The Eilenberg--Steenrod axioms involve a category of pairs.  For a
closed model category $\CAT$, an appropriate category of pairs is the
category of arrows in $\CAT$: A pair is a map $A\to X$ in $\CAT$ and
a map of pairs is a commutative square.  Standard notation is to write
$(X,A)$ for the pair $f\colon A\to X$ (with $f$ understood) and to
write $X$ rather than $(X,\initial)$ for the pair $\initial\to X$
where $\initial$ is the initial object of $\CAT$.  A map of pairs
$(X,A)\to (Y,B)$ is called a weak equivalence when the maps $A\to B$
and $X\to Y$ are weak equivalences in $\CAT$.  In this terminology, we
understand cohomology theories as follows:

\begin{intdefn}
Let $\CAT$ be a closed model category.
A cohomology theory on $\CAT$
consists of a contravariant functor $h^{*}$ from the category of pairs to the
category of graded abelian groups together with natural
transformations of abelian groups $\delta^{n} \colon h^{n}(A)\to
h^{n+1}(X,A)$ for all $n$, satisfying the following axioms:
\begin{enumerate}
\item (Homotopy) If $(X,A)\to (Y,B)$ is a weak equivalence of pairs,
then the induced map $h^{*}(Y,B)\to h^{*}(X,A)$ is an isomorphism of
graded abelian groups.
\item (Exactness) For any pair $(X,A)$, the sequence
\[
\cdots \to h^{n}(X,A)\to h^{n}(X)\to h^{n}(A)\overto{\delta^{n}} h^{n+1}(X,A)
\to \cdots
\]
is exact.
\item (Excision) If $A$ is cofibrant, $A\to B$ and $A\to X$ are
cofibrations, and $Y$ is the pushout $X\cup_{A}B$, then the map of
pairs $(X,A)\to (Y,B)$ induces an isomorphism of graded abelian groups
$h^{*}(Y,B)\to h^{*}(X,A)$.
\item \label{productaxiom}(Product) If $\{X_{\alpha}\}$ is a set of
cofibrant objects and $X$ is the coproduct, then the natural map
$h^{*}(X)\to \prod h^{*}(X_{\alpha})$ is an isomorphism. 
\end{enumerate}
A map of cohomology theories $\phi \colon h^{*}\to k^{*}$ is a natural transformation of
contravariant functors that makes the diagram
\[
\xymatrix@R-=3ex{
h^{n}(A)\ar[r]^{\delta^{n}}\ar[d]_{\phi_{A}}
&h^{n+1}(X,A)\ar[d]^{\phi_{(X,A)}}\\
k^{n}(A)\ar[r]_{\delta^{n}}
&k^{n+1}(X,A)
}
\]
commute for all $(X,A)$, $n$.

A homology theory consists of a covariant functor $h_{*}$ together with
natural transformations $\partial_{n}\colon h_{n+1}(X,A)\to h_{n}(A)$
satisfying completely analogous axioms, with the Product Axiom
replaced by the following Direct Sum Axiom:
\begin{itemize}
\item [$\text{(\ref{productaxiom})}'$](Direct Sum)
If $\{X_{\alpha}\}$ is a set of
cofibrant objects and $X$ is the coproduct, then the natural map 
$\bigoplus h_{*}(X_{\alpha})
\to h_{*}(X)$ is an isomorphism.
\end{itemize}
\end{intdefn}

As a matter of pure algebra, the axioms above imply the stronger
property that for $A\to B\to X$, the sequence 
\[
\cdots \to h^{n}(X,B)\to h^{n}(X,A)\to h^{n}(B,A)\to h^{n+1}(X,B)
\to \cdots
\]
is exact, where the last map is the composite $h^{n}(B,A)\to
h^{n}(B)\to h^{n+1}(X,B)$.  Likewise, they imply the stronger
property that the natural map $h^{*}(\coprod X_{\alpha}, \coprod
A_{\alpha})\to \prod h^{*}(X_{\alpha},A_{\alpha})$ is an isomorphism
when each $X_{\alpha}$ and $A_{\alpha}$ is cofibrant.  
As a matter of pure homotopy theory, whether or not
$(h^{*},\delta^{*})$ or $(h_{*},\partial_{*})$ satisfies the axioms
for a particular model structure on $\CAT$ depends only on the weak
equivalences and not on the cofibrations (see
Section~\ref{secpairs} for details).

\medskip

The purpose of this paper is to study homology and cohomology theories
on categories of \einf ring spectra, or equivalently, on the modern
categories of EKMM commutative $S$-algebras \cite{ekmm}, where the
initial object is the sphere spectrum $S$ and the coproduct is the
modern symmetric monoidal smash product.  Because the final object in
the category of commutative $S$-algebras is the trivial (one-point)
spectrum $*$, if we take $\CAT$ above to be the model category of
commutative $S$-algebras, then there are no non-trivial cohomology
theories: If $A$ is any cofibrant commutative $S$-algebra and $C\to *$ is a
cofibrant approximation, then $A\sma C$ is
contractible, and so by the Homotopy Axiom and the Exactness Axiom,
$h^{*}(A\sma C,C)=0$, for any cohomology theory $h^{*}$.  Then the
Excision Axiom applied to the cofibrations $S\to A$ and $S\to C$
implies that the map
\[
h^{*}(A)=h^{*}(A,S) \to h^{*}(A\sma C,C) =0
\]
is an isomorphism.  Now it easily follows from the Homotopy Axiom and
the Exactness Axiom that $h^{*}$ is zero on any pair.
In order to have non-trivial homology and cohomology theories, we
therefore need to
consider categories of commutative $S$-algebras with non-trivial final
objects.  We do this by considering over-categories, and we work more
generally with categories of commutative $R$-algebras for a
cofibrant commutative $S$-algebra $R$.
For a commutative $R$-algebra $B$, let $\Com{R}{B}$ denote the
category of commutative $R$-algebras lying over $B$.
This is a closed model category with initial object $R$, final object
$B$, and coproduct $\sma_{R}$, the smash product over $R$.

Topological Andr\'e--Quillen cohomology with various coefficients
provides examples of cohomology theories on $\Com{R}{B}$.  For a
cofibrant commutative $R$-algebra $A$ and a cofibrant commutative
$A$-algebra $X$, Basterra \cite{mbthesis} constructs the cotangent
complex $\cotan{A}{X}$ as the derived commutative $X$-algebra
indecomposables (of the derived augmentation ideal) of $X\sma_{A}X$.
The cotangent complex $\cotan{A}{X}$ is an $X$-module, and restricting
to maps $A\to X$ in $\Com{R}{B}$, we can regard
\[ \Bcotan{A}{X}=B\sma_{X}\cotan{A}{X} \]
as a functor from the homotopy category of $\Com{R}{B}$ to the
homotopy category of $B$-modules.  Topological Andr\'e--Quillen
cohomology with coefficients in a $B$-module $M$ is defined by
\[
\TAQ_{R}^{*}(X,A;M) = \Ext^{*}_{B}(\Bcotan{A}{X},M)
\]
(where $\Ext$ is as in \cite[IV.1.1]{ekmm}).
Using the connecting homomorphism in the transitivity sequence
\cite[4.4]{mbthesis}, $\TAQ_{R}^{*}(-;M)$ becomes a cohomology theory on
$\Com{R}{B}$, functorially in $M$ in the homotopy category of
$B$-modules.  Our main result is the following theorem, which says in
particular that every cohomology theory on $\Com{R}{B}$ is isomorphic
to $\TAQ_{R}^{*}(-;M)$ for some $B$-module $M$.

\begin{intthm}\label{intmainthm}
Topological Andr\'e--Quillen cohomology, viewed as a functor from the
homotopy category of $B$-modules to the category of cohomology
theories on $\Com{R}{B}$, is an equivalence of categories.
\end{intthm}

A characterization of the category of homology theories on
$\Com{R}{B}$ is slightly trickier because of the failure of
Brown's Representability Theorem for homology theories
\cite{brownfail}.  Given a homology theory $h_{*}$ on the category
$\Mod{B}$ of $B$-modules, we obtain a homology theory $\newh{h}{*}$ on
the category $\Com{R}{B}$ by setting
\[
\newh{h}{*}(X,A) = h_{*}(\Bcotan{A}{X},*).
\]
This describes a functor from the category of homology theories on
$\Mod{B}$ to the category of homology theories on $\Com{R}{B}$.  We
prove that this functor is an equivalence of categories.

\begin{intthm}\label{inthomthm}
The category of homology theories on $\Com{R}{B}$ is equivalent to the
category of homology theories on $\Mod{B}$.
\end{intthm}

As always, there is a close relationship between homology and
cohomology theories and a category of ``spectra'' that we review in
Section~\ref{secreduced}.  For $R=B$, the category $\Com{B}{B}$ is
enriched over the category of based spaces, with tensors and
cotensors, and so in particular we have a suspension functor.  We
denote this suspension functor by $\susp{B}$ to avoid confusion with
suspension of the underlying $B$-module.  A $\Com{B}{B}$-spectrum is a
sequence of objects $A_{n}$, $n\geq 0$, together with ``structure
maps'' $\sigma_{n}\colon \susp{B}A_{n}\to A_{n+1}$.  A map of spectra
is a collection of maps $A_{n}\to A'_{n}$ that commute with the
structure maps.  We define the homotopy groups of a spectrum
$\csp{A}=\{A_{n}\}$ by
\[
\pi_{q}\csp{A}=\Colim \tpi_{q+n}A_{n},
\]
where $\tpi_{*} A=\Ker(\pi_{*}A\to \pi_{*}B)$.
Standard techniques \cite{hoveyssmc,mmss} allow us to prove in
Section~\ref{secpfspectra} that the category of $\Com{B}{B}$-spectra
forms a closed model category, with weak equivalences the maps that
induce isomorphisms on homotopy groups.  The resulting homotopy
category is called the ``stable category'' of $\Com{B}{B}$.  In
Section~\ref{secspectra}, we prove the following theorem:

\begin{intthm}\label{intstablethm}
Let $B$ be a cofibrant commutative $R$-algebra.  The stable category
of $\Com{B}{B}$ is equivalent to the homotopy category of $B$-modules.
\end{intthm}

In fact, as explained in Section~\ref{secspectra}, the equivalence of
homotopy categories arises from a Quillen equivalence.  Although the
technical hypotheses do not quite apply, this theorem is closely
related to the title theorem of Schwede--Shipley \cite{ssstable} that
stable categories are categories of modules.
Theorem~\ref{intstablethm} is in marked contrast to the corresponding
situation for simplicial commutative algebras studied by Schwede
\cite{sstalgth}, where the stable category is equivalent to the
homotopy category of modules over a ring spectrum that is generally
very different from the ground ring; see also
Theorem~\ref{thmfreestab} below.

For an object $A$ in $\Com{R}{B}$, $B\sma_{R}A$ is naturally an object
of $\Com{B}{B}$ and we have an associated $\Com{B}{B}$-spectrum
$\cspsrel{A}=\{\susp{B}^{n}(B\sma_{R}A)\}$ called the ``suspension
spectrum''.  Closely related to the previous theorems is the following:

\begin{intthm}\label{intcotanc}
Let $B$ be a cofibrant
commutative $R$-algebra, and $A$ a cofibrant object in $\Com{R}{B}$.
Under the equivalence of Theorem~\ref{intstablethm}, the suspension
spectrum $\cspsrel{A}$ corresponds to the $B$-module $\Bcotan{R}{A}$.
\end{intthm}

The suspension for commutative $S$-algebras turns out to be closely
related to delooping for spaces.  When $X$ is an \einf space, that is,
a space with an action of an \einf operad, the work of May, Quinn, and
Ray \cite[IV\S1]{mayeinf} implies that the suspension spectrum
$\Sigma^{\infty}X_{+}$ is an \einf ring spectrum.  We can make sense
of the Andr\'e--Quillen Cohomology and the cotangent complex of \einf
ring spectra. Up to equivalence, these do 
not depend on the operad, and we can understand these in terms of an
equivalent commutative $S$-algebra.  The work of May and Thomason
\cite{maythomason} shows that up to isomorphism in the stable
category, there is a canonical spectrum associated to $X$ whose zeroth
space is the group completion of $X$; it is any spectrum output by an
``infinite loop space machine''. In Section~\ref{secdeloop}, we
reinterpret part of the proof of the previous theorem to prove the
following delooping theorem:

\begin{intthm}\label{intdeloopthm}
For an \einf space $X$, the
cotangent complex of the \einf ring spectrum $\Sigma^{\infty}X_{+}$ is
the extended $\Sigma^{\infty}X_{+}$-module $(\Sigma^{\infty}X_{+})\sma Z$,
where $Z$ is the spectrum associated to $X$.
\end{intthm}

According to Lewis \cite[\S IX]{lms}, the Thom spectrum $M$ obtained from
a map of \einf spaces $X\to BF$ naturally has the structure of an
\einf ring spectrum (see also Mahowald \cite{mahowaldthom}), and the
diagonal map  
\[
M\to M\sma X_{+}
\]
is a map of \einf ring spectra.  The derived extension of scalars to
\einf $M$-algebras,
\[ M\sma M\to M\sma X_{+} \iso M \sma \Sigma^{\infty}X_{+} \]
induces the Thom isomorphism and is a weak equivalence.  The
cotangent complex commutes with extension of scalars
\cite[4.5]{mbthesis}, and we obtain the following corollary of the
previous theorem:

\begin{intcor}\label{intthomcor}
Let $\alpha \colon X\to BF$ be a map of \einf spaces, $M$ the
associated \einf ring Thom spectrum, and $Z$ the spectrum associated
to $X$.  Then the cotangent complex $\cotan{S}M$
is the $M$-module $M\sma Z$.
\end{intcor}

As a special case, we obtain the following result.  In it $bu$ denotes
the spectrum with zeroth \einf space $BU$; it is equivalent to
$\Sigma^{2}ku$, where $ku$ is connective $K$-theory.

\begin{intcor}
The cotangent complex 
$\cotan{S}MU$ of $MU$ is the $MU$-module $MU\sma bu$.
\end{intcor}

We have stated these results in terms of EKMM commutative
$S$-algebras, but because of the homotopy invariant nature of
Theorems~\ref{intmainthm} and~\ref{inthomthm}, they hold in the
context of commutative symmetric spectra and commutative orthogonal
spectra \cite{mmss}, or in any Quillen equivalent modern category of
\einf ring spectra.

More generally, we have stated results in terms of commutative
$S$-algebras, but results analogous to
Theorems~\ref{intmainthm}--\ref{intcotanc} hold for any sort of
operadic algebras in EKMM $S$-modules with slight modifications of the
arguments below; see Section~\ref{secoperadic} for details.

\section{Reduced Theories}\label{secreduced}

Although the most natural statement of Theorems~\ref{intmainthm}
and~\ref{inthomthm} is in terms of cohomology and homology theories
defined on pairs, the most natural proof is in terms of
reduced cohomology and homology theories.  The purpose of this section
is to record some basic facts about the relationship of cohomology
theories on pairs, reduced cohomology theories, ``Omega weak spectra''
(see Definition~\ref{defowksp} below), and spectra.  These
results hold quite generally and their arguments depend very little on
the specifics of the category $\Com{R}{B}$ of commutative $R$-algebras
over $B$.  Since the arguments are familiar from other contexts, we
omit many of the details.

We begin with the definition of reduced theories.  These are defined
for the homotopy categories of ``pointed closed model categories''
(closed model categories where the initial object is the final
object); these categories have the extra structure of a ``suspension
functor'' \cite[I\S2]{quil} and of ``cofiber sequences''
\cite[I\S3]{quil}. 

\begin{defn}\label{defreduced}
Let $\CAT$ be a pointed closed model category.  A reduced cohomology
theory on $\CAT$ consists of a contravariant functor $\rc{h}{*}$ from
the homotopy category $\h\CAT$ to the category of graded abelian
groups together with a natural isomorphisms of abelian groups
$\sic\colon \rc{h}{n}(X)\to \rc{h}{n+1}(\Sigma X)$ (the suspension
isomorphism) for all $n$, satisfying the following axioms:
\begin{enumerate}
\item (Exactness) If $X\to Y \to Z$ is part of a cofibration sequence,
then 
\[
\rc{h}{n}(Z)\to \rc{h}{n}(Y)\to \rc{h}{n}(X)
\]
is exact for all $n$.
\item \label{redproductaxiom}(Product) 
If $\{X_{\alpha}\}$ is a set of objects and $X$ is the coproduct in
$\h\CAT$, then the natural map $\rc{h}{*}(X)\to \prod
h^{*}(X_{\alpha})$ is an isomorphism. 
\end{enumerate}
A reduced homology theory consists of a covariant functor $\rh{h}{*}$
together with natural isomorphisms $\sih\colon \rh{h}{n+1}(\Sigma
X)\to \rh{h}{n}(X)$ satisfying an analogous exactness axiom and the
following Direct Sum Axiom:
\begin{itemize}
\item [$\text{(\ref{redproductaxiom})}'$](Direct Sum) If
$\{X_{\alpha}\}$ is a set of objects and $X$ is the coproduct in
$\h\CAT$, then the natural map $\bigoplus \rh{h}{*}(X_{\alpha}) \to
\rh{h}{*}(X)$ is an isomorphism.
\end{itemize}
A map of reduced cohomology theories or of reduced homology theories
is a natural transformation that commutes with the suspension isomorphisms.
\end{defn}

Whenever the final object in a closed model category is cofibrant,
there is a close relationship between cohomology theories and reduced
cohomology theories on the under-category of the final object.
Writing $B$ for the final object and $\CAT \bs B$ for the
under-category of $B$ (for $\CAT=\Com{R}{B}$, the under-category
$\CAT\bs B$ is $\Com{B}{B}$), a cohomology theory $h$ on $\CAT$ leads
to a reduced cohomology theory $\tilde{h}^{*}$ on $\CAT\bs B$ with
$\tilde{h}^{*}(X)=h^{*}(X,B)$ and the suspension isomorphism (for $X$
cofibrant) obtained from the connecting homomorphism $h^{n}(X,B)\to
h^{n+1}(CX,X)$ and the inverse of the excision isomorphism
$h^{n+1}(\Sigma X,B)\to h^{n+1}(CX,X)$, where $X\to CX$ is a Quillen
cone.  Conversely, given a reduced cohomology theory $\tilde{h}^{*}$
on $\CAT \bs B$, we obtain a cohomology theory on $\CAT$ by setting
$h^{n}(X,A)$ to be $\tilde{h}^{*}$ of the homotopy pushout
$B\cup_{A}X$. In general, we have the following proposition:

\begin{prop}\label{propreduced}
Let $\CAT$ be a closed model category with final object $B$,
and let $B'\to B$ be a cofibrant approximation (an acyclic fibration
with $B'$ cofibrant).  The following categories are
equivalent:
\begin{enumerate}
\item The category of cohomology theories on $\CAT$.
\item The category of cohomology theories on $\CAT / B'$.
\item The category of reduced cohomology theories on $(\CAT /B')\bs B'$.
\end{enumerate}
The analogous result holds for homology theories.
\end{prop}

The homotopy category $\hCom{B}$ of $\Com{B}{B}$ together with a
skeleton of the full
subcategory of finite cell commutative $B$-algebras over $B$ satisfy the
hypotheses of Brown~\cite[\S2]{brown} for a ``homotopy category''.
Brown's Abstract Representability Theorem \cite[2.8]{brown} therefore
applies to show that certain functors are representable.  The
following is the relevant special case; a standard homological algebra
argument (the Barratt--Whitehead ``ladder'' argument) reduces its
hypotheses to those of the Representability Theorem.

\begin{prop}\label{propbrownrep}
Let $h$ be a contravariant functor from $\hCom{B}$ to abelian
groups that satisfies hypotheses (i) and (ii) in the definition of
reduced cohomology theory.  Then there exists an object $X_{h}$ in
$\hCom{B}$ and a natural isomorphism of functors $h(-)\iso
\hCom{B}(-,X_{h})$. 
\end{prop}

It follows that when $\rc{h}{*}$ is a reduced cohomology theory on
$\Com{B}{B}$, each functor $\rc{h}{n}$ is representable by an object
$X_{\rc{h}{n}}$. If we write $\dsusp{B}$ for the left derived functor
of $\susp{B}$, the suspension functor on $\hCom{B}$, and
$\dloops{}$ for its right adjoint, the loop functor on
$\hCom{B}$, then the suspension isomorphism 
\begin{multline*}
\hCom{B}(-,X_{\rc{h}{n}})\iso h^{n}(-) \to h^{n+1}(\dsusp{B}-)
   \iso \hCom{B}(\dsusp{B}-,X_{\rc{h}{n+1}})\\
   \iso \hCom{B}(-,\dloops{}X_{\rc{h}{n+1}})
\end{multline*}
induces (by the Yoneda Lemma) an isomorphism (in $\hCom{B}$),
\[
X_{\rc{h}{n}}\to \dloops{}X_{\rc{h}{n+1}}.
\]
This leads to the following definition:

\begin{defn}\label{defowksp}
Let $\CAT$ be a pointed closed model category.  An Omega weak spectrum
in $\CAT$ consists of objects $X_{0}, X_{1},\dots $, and isomorphisms
$\tsigma_{n}\colon X_{n}\to \Omega X_{n+1}$ in $\h\CAT$, where
$\Omega$ denotes the (Quillen) loop functor on $\h\CAT$.  A map of Omega
weak spectra from $\csp{X}=\{X_{n}\}$ to $\csp{Y}=\{Y_{n}\}$ consists
of maps $X_{n}\to Y_{n}$ in $\h\CAT$ that commute with the structure
maps $\tsigma_{n}$.
\end{defn}

The rule $h^{n}(-;\csp{X})=\hCom{B}(-,X_{n})$ defines a functor from the
category of Omega weak spectra to the category of reduced cohomology
theories on $\Com{B}{B}$.  The Yoneda Lemma implies that this functor
is a full embedding, and the discussion above implies that every
reduced cohomology theory is isomorphic to one associated to an Omega
weak spectrum.  In summary, we have the following proposition:

\begin{prop}\label{propcohowksp}
The category of reduced cohomology theories on $\Com{B}{B}$ is
equivalent to the category of Omega weak spectra in $\Com{B}{B}$.
\end{prop}

Finally, we need a general result about the relationship between Omega
weak spectra and spectra in a simplicial or topological pointed model
category $\CAT$.  In this context, when $X$ is cofibrant, the tensor
$X\sma I_{+}$ is a cylinder object, and so the tensor $\susp{}X=X\sma
S^{1}$ represents the (Quillen) suspension functor on the homotopy
category.  Likewise, when $Y$ is fibrant, the cotensor of $Y$ with
$S^{1}$, $\loops{}Y$, represents the (Quillen) loop functor on the
homotopy category.  A spectrum $\csp{X}$ is defined as a sequence of
objects $X_{0},X_{1},\dots$ and maps $\sigma \colon \susp{}X_{n}\to
X_{n+1}$ in $\CAT$.  We say that a spectrum $\csp{X}$ is ``cofibrant'' if
each $X_{n}$ is cofibrant and each structure map $\susp{}X_{n}\to
X_{n+1}$ is a cofibration.  In the case of interest the cofibrant
spectra in $\Com{B}{B}$ are the cofibrant objects in the stable model
structure (see Theorem~\ref{thmstable} below); quite generally, the
cofibrant spectra are the cofibrant objects in some model structure on
the category of spectra, q.v. \cite[1.13--14]{hoveyssmc}. 

We write $\tsigma$ for the adjoint of the structure map $X_{n}\to
\loops{}X_{n+1}$.  A spectrum $\csp{X}$ is called an Omega spectrum
when each $X_{n}$ is fibrant and each adjoint structure map is a
weak equivalence.  Then by neglect of structure (passing from
$\CAT$ to $\h\CAT$), an Omega spectrum becomes an Omega weak spectrum.
The following lemma is the standard observation that every map
of Omega weak spectra can be rectified to a map of spectra (though
typically not uniquely).

\begin{lem}\label{lemspectrawkspectra}
Let $\CAT$ be a simplicial or topological pointed closed model category.
\begin{enumerate}
\item Every Omega weak spectrum is (non-canonically)
isomorphic to the underlying Omega weak spectrum of a cofibrant Omega
spectrum. 
\item Let $\csp{X}$ and $\csp{Y}$ be Omega spectra and suppose that
$\csp{X}$ is cofibrant.  Any map of Omega weak spectra $f\colon
\csp{X}\to\csp{Y}$ is represented by a map of spectra (generally not
uniquely).
\end{enumerate}
\end{lem}

\begin{proof}
Given an arbitrary Omega weak spectrum $\csp{X}$, each $X_{n}$ is
isomorphic in $\h\CAT$ to an object $X_{n}'$ that is 
fibrant, and using these isomorphisms, we get an isomorphic Omega weak
spectrum $\csp{X}'$.  We choose a cofibrant Omega spectrum $\csp{X}''$
as follows:  We choose $X_{0}''$ to be a cofibrant approximation of
$X_{0}'$.  Then since $X_{0}''$ is cofibrant and $\loops{}X_{1}'$ is
fibrant, we can choose a map $X_{0}''\to \loops{}X_{1}'$ representing
the composite map in $\h\CAT$
\[
X_{0}'' \to X_{0}' \to \Omega X_{1}'.
\]
We factor the adjoint map $\susp{}X_{0}''\to X_{1}'$ as a cofibration
followed by an acyclic fibration
\[
\susp{}X_{0}'' \to X_{1}'' \to X_{1}'.
\]
Continuing in this fashion constructs a cofibrant Omega spectrum
$\csp{X}''$ and an isomorphism of Omega weak spectra $\csp{X}''\to
\csp{X}'$. 

For (ii), since $X_{0}$ is cofibrant and $Y_{0}$ is fibrant, we can
choose a map $X_{0}\to Y_{0}$ representing $f_{0}$.  The
hypothesis SM7 that $\CAT$ is
a simplicial model category or the analogous hypothesis that $\CAT$ is
a topological model category implies that the induced
map of simplicial sets or of spaces
\[
\CAT(X_{1},Y_{1})\to \CAT(\susp{}X_{0},Y_{1})
\]
is a fibration and identifies the induced map on components as 
$\h\CAT(X_{1},Y_{1})\to \h\CAT(\susp{}X_{0},Y_{1})$.  Since
the given map $f_{1}$ in $\h\CAT$ maps to the same component as the
map adjoint to the composite $X_{0}\to Y_{0}\to \loops{}Y_{1}$, there
exists a map $X_{1}\to Y_{1}$ in $\CAT$ that represents $f_{1}$ and
makes the diagram
\[ \xymatrix{
\susp{}X_{0}\ar[r]\ar[dr]
&X_{1}\ar[d]\\
&Y_{1}
} \]
commute in $\CAT$.  Applying this argument inductively to each map
$X_{n}\to Y_{n}$ constructs a compatible map $X_{n+1}\to Y_{n+1}$, and
the map of spectra $\csp{X}\to \csp{Y}$.
\end{proof}

\section{Nucas, Indecomposables, and Stabilization}\label{secnuca}

The functors $\cotan{A}{X}$ and $\Bcotan{A}{X}$ mentioned in the
introduction are somewhat awkward to work with formally because they
are composites of both left and right derived functors.  However,
the technical trick introduced in \cite{mbthesis} of working
with non-unital commutative algebras (``nucas'') alleviates this
problem by providing a point-set (left adjoint) functor whose left
derived functor is a model for $\Bcotan{}$.  The first half of
this section consists of a brief overview of the theory of nucas from
\cite{mbthesis}. Another technical
advantage of the category of nucas is that every object is fibrant,
and this makes the construction of a ``stabilization'' functor $\stab$
from cofibrant spectra to Omega spectra easier.  The second half of
this section constructs this functor and explores its basic properties.

\begin{defn}
Let $B$ be a commutative $R$-algebra.
A non-unital commutative $B$-algebra (or $B$-nuca) consists of a
$B$-module $N$ together with an associative and commutative
multiplication $\mu\colon N\sma_{B}N\to N$.    A map of $B$-nucas
is a map of $B$-modules $N\to N'$ that commutes with the
multiplications.  We write $\NUCA{B}$ for the category 
of $B$-nucas.
\end{defn}

The category of $B$-nucas may also be described as the category of
algebras in $B$-modules over the operad $\oNCom$ with $\oNCom(k)=*$
for $k>0$ and $\oNCom(0)$ empty.  We have a free $B$-nuca functor
\[
\bNCom M = \bigvee_{k>0} M^{(k)}/\Sigma_{k}
\]
for a $B$-module $M$ (where $M^{(k)}=M\sma_{B}\cdots \sma_{B}M$), and
we have a functor $K$ from 
$\NUCA{B}$ to $\Com{B}{B}$ defined by formally adding a unit: On the
underlying $B$-modules, 
\[
KN = B \vee N,
\]
with the multiplication extended from $N$ to $KN$ by the
usual multiplication on $B$ and the $B$-action maps of $B$ on $N$.
The functor $K$ is the left adjoint of the ``augmentation ideal''
functor $I$ from $\Com{B}{B}$ to $\NUCA{B}$ defined by setting $I(A)$
to be the (point-set) fiber of the augmentation $A\to B$. According to
Basterra \cite[1.1,2.2]{mbthesis}, the adjunction $(K,I)$ is a Quillen
equivalence. 

\begin{prop}\label{propki}
The category $\NUCA{B}$ is a topological pointed closed model
category with weak equivalences and fibrations the maps that are weak
equivalences and fibrations (resp.) of the underlying $B$-modules.
The adjunction $(K,I)$ is a Quillen equivalence between $\NUCA{B}$ and
$\Com{B}{B}$. 
\end{prop}

Since $K$ preserves all weak equivalences, it is harmless
to use the same notation for the derived functor.  We denote the right
derived functor of $I$ by $\derI$.  As in any Quillen adjunction,
the derived functor $\derI$ preserves the Quillen loop functor, but
since the derived functor $K$ is an inverse equivalence, it also
preserves the Quillen loop functor, and we obtain the following
proposition:

\begin{prop}\label{propowe}
The derived functors $K$ and $\derI$ induce an equivalence between the
category of Omega weak spectra in $\NUCA{B}$
and the category of Omega weak spectra in $\Com{B}{B}$.
\end{prop}

Basterra \cite[\S3]{mbthesis} constructs an ``indecomposables''
functor $\indec{B}$ from $\NUCA{B}$ to $\Mod{B}$ that is left adjoint to
the ``zero multiplication'' functor $\zero{B}$ from $\Mod{B}$ to
$\NUCA{B}$.  Precisely, the indecomposables functor $\indec{B}N$ is
defined as the (point-set) pushout 
\[
\indec{B}N = *\cup_{(N\sma_{B}N)}N
\]
in the category of $B$-modules of the multiplication $N\sma_{B}N\to N$
over the trivial map $N\sma_{B}N\to *$.  The functor $\zero{B}$ simply assigns
a $B$-module $M$ the trivial map $M\sma_{B}M\to M$ for its multiplication.
Since $\zero{B}$ clearly preserves weak equivalences and fibrations,
we have the following observation of \cite{mbthesis}:

\begin{prop}\label{propqz}
The functors $(\indec{B},\zero{B})$ between $\NUCA{B}$ and $\Mod{B}$
form a Quillen adjunction.
\end{prop}

Again, since $\zero{B}$ preserves all weak equivalences, it is harmless to
denote its derived functor by the same notation.  We write
$\derindec{B}$ for the left derived functor of $\indec{B}$.

The derived functors $\derindec{}$ and $\derI$ are needed in
\cite[4.1]{mbthesis} to defined the cotangent complex.  For a
cofibrant commutative $R$-algebra $A$ and a cofibration of commutative
$R$-algebras $A\to X$, the cotangent complex of $X$ relative to $A$ is
defined to be the $X$-module
\[
\cotan{A}{X} = \derindec{X}\derI (X\sma_{A} X),
\]
where $\derI$ and $\derindec{X}$ are understood in terms of
$\Com{X}{X}$ and $\NUCA{X}$.  When $A$ is not cofibrant or $A\to X$ is
not a cofibration, the cotangent complex is constructed by choosing a
cofibrant approximation $A'\to X'$ of $A\to X$, and then using the
derived extension of scalars functor $X\sma_{X'}^{\dL}(-)$:
\[
\cotan{A}{X} = X\sma_{X'}^{\dL} \cotan{A'}{X'}
= X\sma_{X'}^{\dL}(\derindec{X'}\derI (X'\sma_{A'} X')).
\]
Because the category this construction lands in depends on its inputs,
a discussion of functoriality would be somewhat involved, but (as
mentioned in the introduction), for $(X,A)$ a pair in $\Com{R}{B}$,
the construction we consider,
\[
\Bcotan{A}{X} = B \sma^{\dL}_{X'} \cotan{A'}{X'}
\]
assembles to a functor from the category of pairs in $\Com{R}{B}$ to
the homotopy category of $B$-modules, and this functoriality suffices
for our work.   The last fact we need about the cotangent complex is
the following version of \cite[4.4--5]{mbthesis}:

\begin{prop}\label{propbcotan}
Assume that $B$ is a cofibrant commutative $R$-algebra.  If $A$ is
cofibrant and $A\to X$ is a cofibration in $\Com{R}{B}$, 
then $\Bcotan{A}{X}$ is isomorphic to $\derindec{B}\derI
(B\sma_{A}X)$, naturally in $A$ and cofibrations $A\to X$.
\end{prop}

Next we move on to the stabilization functor.  In other contexts, this
functor is typically denoted by $Q$, but here we denote it by $\stab$
to avoid confusion with the indecomposables functor.  We have the
notion of a spectrum in the category of $B$-nucas, as discussed in the
previous section.  The stabilization functor $\stab$ turns out to be a
functor from spectra in $\NUCA{B}$ to Omega spectra in $\NUCA{B}$.

\begin{defn}\label{defstab}
For a $\NUCA{B}$-spectrum $\csp{X}$, we define a $\NUCA{B}$-spectrum 
$\stab\csp{X}=\{\nstab_{n}\csp{X}\}$ as follows:
Let
\[
\nstab_{n}\csp{X} = \Tel_{k\geq n} \loops{}^{k-n}X_{k},
\]
the telescope over the adjoint structure maps.
We have a map of telescopes 
\[
\Tel_{k\geq n+1} \susp{B}\loops{}^{k-n}X_{k} \to 
\Tel_{k\geq n+1}  \loops{}^{k-(n+1)}X_{k}
\]
induced by sending $\susp{B}\loops{}^{k-n}X_{k}$ to
$\loops{}^{k-(n+1)} X_{k}$ using the counit of the suspension, loop 
adjunction applied to the innermost factor of $\loops{}$, and we have
a map of telescopes
\[
\Tel_{k\geq n} \susp{B}\loops{}^{k-n}X_{k} \to 
\Tel_{k\geq n+1} \susp{B}\loops{}^{k-n}X_{k} 
\]
induced by collapsing down the map $\susp{B}X_{n}\to \susp{B}\loops{}
X_{n+1}$.  We define the structure map $\sigma \colon
\susp{B}\nstab_{n}\csp{X}\to \nstab_{n+1}\csp{X}$ to be the composite
\begin{multline*}
 \susp{B}(\Tel_{k\geq n}\loops{}^{k-n}X_{k}) \iso
\Tel_{k\geq n} \susp{B}\loops{}^{k-n}X_{k} \to \\
\Tel_{k\geq n+1} \susp{B}\loops{}^{k-n}X_{k} \to 
\Tel_{k\geq n+1}  \loops{}^{k-(n+1)}X_{k}.
\end{multline*}
\end{defn}

The construction $\stab$ assembles in the obvious way to a functor
from $\NUCA{B}$-spectra to $\NUCA{B}$-spectra.   
Since the composite map 
\[
\susp{B}X_{n}\overto{\susp{} \tsigma} \susp{B}\loops{}X_{n+1}\to X_{n+1}
\]
is the structure map $\sigma$, the inclusion of $X_{n}$ into the
telescope defining $\nstab_{n}\csp{X}$ defines a natural
transformation $\csp{X}\to \stab\csp{X}$.  The main fact
we need about $\stab$ is the following proposition.  In it, the
homotopy groups of a spectrum $\csp{X}$ are defined by
\[
\pi_{q}\csp{X}=\Colim \pi_{q+n}X_{n}.
\]

\begin{prop}\label{propstab}
For any spectrum $\csp{X}$ in $\NUCA{B}$, $\stab\csp{X}$ is an Omega
spectrum in $\NUCA{B}$ and the natural transformation $\csp{X}\to
\stab\csp{X}$ induces an isomorphism on homotopy groups.
\end{prop}

\begin{proof}
The telescope in the category of $B$-nucas is naturally weakly
equivalent to the telescope in the category of $B$-modules, and so the
usual map 
\[
\loops{}(\Tel_{k\geq n+1} \loops{}^{k-(n+1)}X_{k})\to
\Tel_{k\geq n+1} \loops{}^{k-n}X_{k}
\]
is a weak equivalence.  The composite of the adjoint structure map and
the map above, 
\[
\Tel_{k\geq n} \loops{}^{k-n}X_{k}\to
\loops{}(\Tel_{k\geq n+1} \loops{}^{k-(n+1)}X_{k})\to
\Tel_{k\geq n+1} \loops{}^{k-n}X_{k}
\]
is a homotopy equivalence, and so the adjoint structure map
$\nstab_{n} \csp{X}\to \loops{}\nstab_{n+1}\csp{X}$ is a weak
equivalence.  Thus $\stab\csp{X}$ is an 
Omega spectrum.  We have 
\[
\pi_{q+n}\nstab_{n}\csp{X} \iso \Colim_{k\geq n} \pi_{q+k} X_{k},
\]
and the map $X_{n}\to \nstab_{n}\csp{X}$ induces on homotopy groups
the inclusion of $\pi_{q+n}X_{n}$ into this colimit system.  Under
this identification, the map $\csp{X}\to \stab\csp{X}$ induces on
homotopy groups the map
\[
\Colim_{n}\pi_{q+n}X_{n} \to 
\Colim_{n}\Colim_{k\geq n} \pi_{q+k} X_{k},
\]
which is clearly an isomorphism.
\end{proof}

For a $B$-nuca $N$, let $\csps{N}$ be the ``suspension spectrum''
which has as its $n$-th object the $n$-fold suspension $\susp{B}^{n}N$
and structure maps the identity map
\[
\susp{B}(\susp{B}^{n}N) \to \susp{B}^{n+1}N.
\]
In the case when $N$ is free, i.e., $N=\bNCom M$ for some $B$-module
$M$, we have a canonical isomorphism 
\[
\susp{B}^{n}\bNCom M \iso \bNCom \Sigma^{n}M,
\]
and in particular, we have canonical maps $\Sigma^{n}M\to
\nstab_{n}\csps{\bNCom M}$ for all $n$.
The following result on the suspension spectra of free
$B$-nucas represents the fundamental difference between the
context of commutative $S$-algebras and simplicial commutative algebras.

\begin{thm}\label{thmfreestab}
Assume that $B$ is a cofibrant commutative $R$-algebra.  If $M$ is a
cofibrant $B$-module, then the canonical maps $\Sigma^{n}M\to
\nstab_{n}\csps{\bNCom M}$ are weak equivalences for all $n$.
\end{thm}

\begin{proof}
The general case follows from the case $n=0$, where we are studying the map 
\[
\pi_{q}\bNCom M\iso
\pi_{q+k}\Sigma^{k}\bNCom M\to
\pi_{q+k}\bNCom \Sigma^{k}M.
\]
The map $\Sigma^{k}\bNCom M\to \bNCom \Sigma^{k}M$ takes the wedge
summand $\Sigma^{k}M^{(m)}/\Sigma_{m}$ to the corresponding wedge
summand $(\Sigma^{k}M)^{(m)}/\Sigma_{m}$ via the diagonal map on the
sphere $S^{k}$.  The proposition is an immediate consequence of the
following lemma.
\end{proof}

\begin{lem}\label{lemsymprod}
Let $B$ be a cofibrant commutative $R$-algebra, let $M$ be a cofibrant
$B$-module, and let $x$ be an element of $\pi_{q}(M^{(m)}/\Sigma_{m})$
for some $m>1$ and some integer $q$.  Then for some $k$, the composite
map 
\[
\pi_{q}(M^{(m)}/\Sigma_{m})\iso
\pi_{q+k}\Sigma^{k}(M^{(m)}/\Sigma_{m})\to
\pi_{q+k}((\Sigma^{k}M)^{(m)}/\Sigma_{m})
\]
sends $x$ to zero.
\end{lem}

\begin{proof}
By \cite[III.5.1]{ekmm}, for any cofibrant $B$-module $N$ (e.g., $M$,
$\Sigma^{k}M$), the map
\[
E\Sigma_{m+}\sma_{\Sigma_{m}}N^{(m)}\to N^{(m)}/\Sigma_{m}
\]
is a weak equivalence.  The cellular filtration of the $\Sigma_{m}$-CW
complex $E\Sigma_{m}$ induces an increasing filtration on
$E\Sigma_{m+}\sma_{\Sigma_{m}}N^{(m)}$ and on the homotopy groups of
$N^{(m)}/\Sigma_{m}$; clearly, this filtration is trivial (zero) below
the zero level.  The map 
\[
\delta \colon E\Sigma_{m+}\sma_{\Sigma_{m}}\Sigma M^{(m)}\to
E\Sigma_{m+}\sma_{\Sigma_{m}}(\Sigma M)^{(m)}
\]
(induced by the diagonal $S^{1}\to (S^{1})^{(m)}$)
preserves the filtration.  Since $m>1$, the map $\Sigma M^{(m)}\to (\Sigma
M)^{m}$ is null homotopic, and so $\delta$ induces the zero map on the
$E^{1}$-term of the homotopy group spectral sequence associated to the
filtration.  It follows that the map on homotopy groups induced by
$\delta$ strictly lowers filtration level.  The map
$\pi_{q}\delta^{k}$ is the map 
\[
\pi_{q}(M^{(m)}/\Sigma_{m})\iso
\pi_{q+k}\Sigma^{k}(M^{(m)}/\Sigma_{m})\to
\pi_{q+k}((\Sigma^{k}M)^{(m)}/\Sigma_{m})
\]
in the statement, which 
therefore takes every element in filtration level $n$ to an element of
filtration level $n-k$.  Taking $k$ greater than the minimum
filtration level of $x$, the map must send $x$ to zero.
\end{proof}

We can illustrate the previous lemma and theorem in the case $R=S$ and
$B=H\bF_{2}$ since in this case, $\pi_{*}\bNCom M$ is easy to describe
in terms of $\pi_{*}M$.  Specifically, $\pi_{*}\bNCom M$ is the
polynomial $\bF_{2}$-nuca on the free allowable Dyer-Lashof (``DL'')
module on $\pi_{*}M$ modulo the relation that the square (in the
algebra structure) is equal to the squaring operation (in the DL
structure) on each element.  The suspension map $\sigma \colon
\pi_{*}\bNCom M\to \pi_{*+1}\bNCom \Sigma M$ kills decomposables (in
the algebra structure) and is a map of DL-modules.  An element of the
form $Q_{s}x$ for $x\in \pi_{q}\bNCom M$ is therefore killed by the
map $\sigma^{k}\colon \pi_{q}\bNCom M\to \pi_{q+k}\bNCom \Sigma^{k}M$
for $k=s-q+1$ because
$\sigma^{k-1}(Q_{s}x)=Q_{s}(\sigma^{k-1}x)=(\sigma^{k-1}x)^{2}$ is
decomposable.  From this it is easy to see directly that the map
$\pi_{*}M\to \Colim \pi_{*+k}\bNCom \Sigma^{k}M$ is an isomorphism.

\section{Proof of Theorems~\ref{intstablethm} and \ref{intcotanc}}
\label{secspectra}  

In this section, we prove Theorems~\ref{intstablethm}
and~\ref{intcotanc} of the introduction.  The arguments take
advantage of the technical simplifications the category of $B$-nucas
provide and use Theorem~\ref{thmfreestab} above as the key step.  They
also make use of the following theorem, proved in
Section~\ref{secpfspectra}:

\begin{thm}\label{thmstable}
Let $\CAT$ be $\Mod{B}$, $\NUCA{B}$, or $\Com{B}{B}$.  Then the
category $\spc{\CAT}$ of $\CAT$-spectra is a topological closed
model category with:
\begin{enumerate}
\item Cofibrations the maps $\csp{X}\to \csp{Y}$ with $X_{0}\to Y_{0}$
a cofibration and each $\susp{}Y_{n}\cup_{\susp{} X_{n}}X_{n+1}\to Y_{n+1}$ a
cofibration,
\item Fibrations the maps $\csp{X}\to \csp{Y}$ with each $X_{n}\to
Y_{n}$ a fibration and each $X_{n}\to
Y_{n}\times_{\loops{}Y_{n+1}}\loops{}X_{n+1}$ 
a weak equivalence, and
\item Weak equivalences the maps that induce an isomorphism on
homotopy groups. 
\end{enumerate}
\end{thm}

In the statement, ``$\susp{}$'' and ``$\loops{}$'' denote the
(point-set) tensor and cotensor with the based space $S^{1}$ in the
pointed topological category $\CAT$, and ``cofibration'' means a
cofibration in the model structure on $\CAT$ (which was called a
``q-cofibration'' in \cite{mbthesis} and \cite{ekmm}).  We have
defined homotopy groups for $\Com{B}{B}$-spectra and
$\NUCA{B}$-spectra 
above, and the definition is the same for $\Mod{B}$-spectra:
\[
\pi_{q}\csp{X}=\Colim \pi_{q+n}X_{n}.
\]
To avoid confusion, we use the term ``stable equivalence'' for weak
equivalence in the model structure on $\spc{\CAT}$ above, and we call
the homotopy category of this model structure the ``stable category''
of $\CAT$.

We have the following easy consequence of the characterization of
fibrations: 

\begin{prop}\label{propom1}
An object is fibrant in one of the model structures in
Theorem~\ref{thmstable} if and only if it is an Omega spectrum. 
\end{prop}

The following proposition is also clear:

\begin{prop}\label{propom2}
A map of Omega spectra $\csp{X}\to \csp{Y}$ is a stable equivalence if
and only if it is a weak equivalence $X_{0}\to Y_{0}$.
\end{prop}

For any of the model categories $\CAT$ in Theorem~\ref{thmstable}, we
have a Quillen adjunction between $\CAT$ and the category $\spc{\CAT}$ of
$\CAT$-spectra with left adjoint the suspension spectrum functor
(that sends an object $X$ to the suspension spectrum
$\sps{X}=\{\susp{}^{n}X\}$) and the zeroth object functor (that sends
a spectrum $\csp{X}$ to the zeroth object $X_{0}$).
The previous propositions applied to $\CAT=\Mod{B}$ imply that this is
a Quillen equivalence.

\begin{prop}\label{propszmod}
Let $\CAT$ be one of the categories in Theorem~\ref{thmstable}.  The
suspension spectrum functor and the zeroth object functor are a
Quillen adjunction between the model category $\CAT$ and the model
category $\spc{\CAT}$ of $\CAT$-spectra.  In the case of
$\CAT=\Mod{B}$, this Quillen adjunction is a Quillen equivalence.
\end{prop}

If we consider any adjunction that is enriched over the category of
based spaces, the left adjoint preserves tensors and the right adjoint
preserves cotensors, and so both functors extend to functors between
the categories of spectra.  It is easy to see that the induced
functors on spectra remain adjoints.  When in addition the categories
are ones considered in Theorem~\ref{thmstable} and the adjunction is a
Quillen adjunction, the characterization of the cofibrations and
fibrations in the categories of spectra imply that the induced
adjunction on spectra is a Quillen adjunction.  

When we combine the previous proposition with the observations of the
last paragraph applied to the free, forgetful adjunction
$\Mod{B}\tofrom\NUCA{B}$ and to the $K,I$ adjunction
$\NUCA{B}\tofrom\Com{B}{B}$, we have the following sequence of Quillen
adjunctions:
\begin{equation}\label{dispadj} 
\xymatrix{
&\Mod{B} \ar@<.75ex>[r]^{\Sigma^{\infty}}
&\spc{\Mod{B}}
\ar@<.75ex>[l]^{(-)_{0}}
\ar@<.75ex>[r]^{\csp\bNCom}
&\spc{\NUCA{B}}
\ar@<.75ex>[l]
\ar@<.75ex>[r]^{\csp K}
&\spc{\Com{B}{B}}
\ar@<.75ex>[l]^{\csp I}
} 
\end{equation}
The first and last of these Quillen adjunctions are Quillen
equivalences, and we prove that 
the middle one is also a Quillen equivalence.  

\begin{lem}\label{lemfreeforg}
Assume that $B$ is a cofibrant commutative $R$-algebra.  Then the free,
forgetful adjunction between $\Mod{B}$-spectra and 
$\NUCA{B}$-spectra is a Quillen equivalence.
\end{lem}

\begin{proof}
Let $M$ be a cofibrant $B$-module and let $\csp{X}$ be a spectrum in
$\NUCA{B}$ that is fibrant in the model structure above.  It suffices to
show that a map $\Sigma^{\infty}M\to \csp{X}$ is a stable equivalence if and
only if the adjoint map $\csps{\bNCom M}\to \csp{X}$ is a weak equivalence.
By Proposition~\ref{propstab} and Propositions~\ref{propom1}
and~\ref{propom2} above, it
suffices to show that $M\to X_{0}$ is a weak equivalence if and only
if $\nstab_{0}\csps{\bNCom M}\to \nstab_{0}\csp{X}$ 
is a weak equivalence.  The diagram
\[ \xymatrix{
M\ar[r]\dto&X_{0}\ar[d]\\
\nstab_{0}\csps{\bNCom M}\ar[r]&\nstab_{0}\csp{X}
} \]
in $\Mod{B}$ commutes, and the result follows from
Theorem~\ref{thmfreestab}.
\end{proof}

As an immediate consequence, we
obtain a Quillen equivalence between the category of $B$-modules and
the category of $\Com{B}{B}$-spectra.  This is sufficient to prove 
Theorem~\ref{intstablethm} as stated, but we would
like to know that the equivalence of homotopy categories is induced by
the functor that sends a $B$-module $M$ to the
$\Com{B}{B}$-spectrum 
\[
\fancyzero{M} = \{K\zero{B}\Sigma^{n}M \}
\]
that represents Andr\'e--Quillen cohomology with coefficients in $M$.
Also, for Theorem~\ref{intcotanc}, we need to know that the
equivalence of homotopy categories takes the suspension spectrum
$\cspsrel{A}$ of a cofibrant object $A$ of $\Com{R}{B}$ to the
$B$-module $\Bcotan{R}{A}$.  Both of these are consequences of the
following theorem.

\begin{thm}\label{thmqzequiv}
Let $B$ be a cofibrant commutative $R$-algebra.  Then
the $\csp{\indec{B}}$, $\csp{\zero{B}}$ adjunction between $\spc{\NUCA{B}}$ and
$\spc{\Mod{B}}$ is a Quillen equivalence.
\end{thm}

\begin{proof}
Since the Quillen adjunction $\indec{B}$, $\zero{B}$
between $\NUCA{B}$ and $\spc{\Mod{B}}$ is enriched over the category
of based spaces, as observed above, we obtain a Quillen adjunction
$\csp{\indec{B}}$, $\csp{\zero{B}}$ between the categories of spectra.
According to MMSS \cite[A.2.(ii)]{mmss}, for this to be a Quillen
equivalence, we just need to show that one of the derived functors is
an equivalence on the homotopy categories.  
If we write $\csp{\derindec{B}}$ for the left
derived functor of  $\csp{\indec{B}}\colon
\spc{\NUCA{B}}\to \spc{\Mod{B}}$ and $\csp{\derbNCom}$ for
the left derived functor of $\csp{\bNCom}\colon \spc{\Mod{B}}\to
\spc{\NUCA{B}}$, then the composite functor $\csp{\derindec{B}}\circ
\csp{\derbNCom}$ is naturally isomorphic to the derived functor of the
composite and so is naturally isomorphic to the identity functor on
$\spc{\Mod{B}}$.  By the previous
lemma, $\csp{\derbNCom}$ is an equivalence, and it follows
that $\csp{\derindec{B}}$ is an equivalence. We conclude that 
$\csp{\indec{B}}$, $\csp{\zero{B}}$ is a Quillen equivalence.
\end{proof}

In the course of the previous argument, we also proved the following result:

\begin{prop}\label{propindecforg}
If $B$ is a cofibrant commutative $R$-algebra, then
the derived functors of $\csp\bNCom$ and $\csp{\zero{B}}$ 
are naturally isomorphic functors
$\h\spc{\Mod{B}}\to\h\spc{\NUCA{B}}$, and
the derived functors of
$\csp{\indec{B}}$ and the forgetful 
functor are naturally isomorphic functors
$\h\spc{\NUCA{B}}\to\h\spc{\Mod{B}}$. 
\end{prop}

We can be more explicit about these natural isomorphisms. The first is
induced by the natural transformation of point-set functors $\bNCom\to
\zero{B}$ that arises from the universal property of the free functor.
When $\csp{M}$ is a cofibrant object in $\spc{\Mod{B}}$, the map
$\csp{\bNCom}\csp{M}\to \csp{\zero{B}}\csp{M}$ is a stable equivalence
because the composite map
\[
\csp{M}\to\csp{\bNCom}\csp{M}\to \csp{\zero{B}}\csp{M}
\]
in $\spc{\Mod{B}}$ is the identity and the map
$\csp{M}\to\csp{\bNCom}\csp{M}$ is a stable equivalence by
Lemma~\ref{lemfreeforg}.  The other natural transformation is induced
by the unit of the $\indec{B}$, $\zero{B}$ adjunction, specifically, the
natural transformation $\Id\to \csp{\indec{B}}$ in $\spc{\Mod{B}}$.
This map is hard to study directly; instead, the argument implicit in
the proof of Theorem~\ref{thmqzequiv} studies the (solid) zigzag
\[ \xymatrix{
\csp{X}'\ar[r]\ar[dr]&
\csp\bNCom \csp{X}' \ar[r]\ar@{.>}[d]
&\csp{\indec{B}} \csp\bNCom \csp{X}'\ar[d]\\
&\csp{X} \ar@{.>}[r]& \csp{\indec{B}}\csp{X} 
} 
\]
where $\csp{X}$ is a cofibrant object in $\spc{\NUCA{B}}$, and
$\csp{X}'\to \csp{X}$ is a cofibrant approximation in $\spc{\Mod{B}}$.
The diagonal solid arrow is therefore a stable equivalence by
assumption, and the composite horizontal map $\csp{X}'\to
\csp{\indec{B}} \csp\bNCom \csp{X}'$ is the isomorphism
(explicit) in the proof of the theorem.  The remaining solid arrows
are also stable equivalences:
The top one $\csp{X}'\to \csp\bNCom \csp{X}'$ is a stable equivalence
by the lemma and the right-hand one is a
stable equivalence because, as a Quillen left adjoint,
$\csp{\indec{B}}$ preserves stable equivalences between cofibrant
objects.  This is the concrete argument that shows that the natural map
$\csp{X}\to\csp{\indec{B}}\csp{X}$ in $\spc{\Mod{B}}$ is a stable
equivalence for $\csp{X}$ cofibrant in $\spc{\NUCA{B}}$.

\begin{proof}[Proof of Theorems~\ref{intstablethm} and~\ref{intcotanc}]
The equivalence of the stable category of $\Com{B}{B}$ and the
homotopy category of $B$-modules we need for
Theorem~\ref{intstablethm} follows from the three Quillen equivalences
in \eqref{dispadj}, or better, by the outer two in \eqref{dispadj} and
the one in Theorem~\ref{thmqzequiv}.  Theorem~\ref{intcotanc} then
follows by Proposition~\ref{propbcotan}.
\end{proof}

\section{Proof of Theorem~\ref{intmainthm}}

In this section, we prove Theorem~\ref{intmainthm} from the
introduction.  By Propositions~\ref{propreduced}, \ref{propcohowksp}, and
\ref{propowe}, it suffices to consider the case when $B$ is a
cofibrant commutative $R$-algebra and show that the functor that sends
a $B$-module $M$ to the Omega weak spectrum in $\NUCA{B}$
\[
\csp{Z}\Sigma^{\infty}M = \{Z \Sigma^{n} M\}
\]
gives an equivalence between the homotopy category of $B$-modules and
the category of Omega weak spectra in $\NUCA{B}$. We write $\owspcn$
for the category of Omega weak spectra in $\NUCA{B}$.

The work of the previous section suggests that we should be able to
use the forgetful functor on the zeroth object as the inverse
equivalence $\owspcn\to \h\Mod{B}$.  The composite functor
$\h\Mod{B}\to \h\Mod{B}$ is the identity.  The composite functor
$\owspcn\to\owspcn$ sends $\csp{X}$ to $\csp{Z}\Sigma^{\infty}X_{0}$;
the proof of Theorem~\ref{intmainthm} is completed by showing that
this functor is naturally isomorphic to the identity.

To construct an isomorphism in $\owspcn$ between $\csp{X}$ and
$\csp{Z}\Sigma^{\infty}X_{0}$, we have to take a detour through the
category of spectra in $\NUCA{B}$.  To avoid confusion, we write
$\uw\csp{A}$ for the underlying Omega weak spectrum of an Omega
spectrum $\csp{A}$.  First we choose a cofibrant
Omega spectrum $\spapp{X}$ and an isomorphism of Omega weak spectra
$\uw\spapp{X}\to \csp{X}$ as in Lemma~\ref{lemspectrawkspectra}.(i).
We write $\zspapp{X}$ for the zeroth object of $\spapp{X}$, and
we choose a cofibrant approximation $\cofapp{X}\to \zspapp{X}$
in the category of $B$-modules.
Then we have the following chain of stable equivalences in
$\spc{\NUCA{B}}$: 
\[
\spapp{X}\from \csp\bNCom \Sigma^{\infty}\cofapp{X} \to
\csp{Z}\Sigma^{\infty}\cofapp{X} \to
\csp{Z}\Sigma^{\infty}\zspapp{X}.
\]
(See the explanation following Proposition~\ref{propindecforg} for a
proof that the first two maps are stable equivalences.)
Using the functor $\stab$, we then have the following chain of stable
equivalences of Omega spectra:
\[
\spapp{X}\to \stab\spapp{X}\from 
\stab(\csp\bNCom \Sigma^{\infty}\cofapp{X}) \to 
\stab(\csp{Z}\Sigma^{\infty}\zspapp{X}) \from 
\csp{Z}\Sigma^{\infty}\zspapp{X}.
\]
Finally, applying $\uw$, we have the following chain of isomorphisms
of Omega weak spectra:
\begin{multline*}
\csp{X} \iso
\uw\spapp{X}\to \uw\stab\spapp{X}\from 
\uw\stab(\csp\bNCom \Sigma^{\infty}\cofapp{X}) \to \\
\uw\stab(\csp{Z}\Sigma^{\infty}\zspapp{X}) \from 
\uw\csp{Z}\Sigma^{\infty}\zspapp{X} \iso
\csp{Z}\Sigma^{\infty}X_{0}.
\end{multline*}
Let $\bigiso{X}\colon \csp{X}\to \csp{Z}\Sigma^{\infty}X_{0}$ be the
composite isomorphism.

The isomorphism $\bigiso{X}$ appears to depend on the choices made
above, and it is far from obvious that $\bigisonat$ is a natural
transformation.  Let $f\colon \csp{X}\to \csp{Y}$ be a map of Omega
weak spectra; we must show that the diagram
\[
\xymatrix@R-=3ex{
\csp{X}\ar[rr]^{\bigiso{X}}\ar[d]_{f}
&&\csp{Z}\Sigma^{\infty}X_{0}\ar[d]^{f_{0}}\\
\csp{Y}\ar[rr]_{\bigiso{Y}}
&&\csp{Z}\Sigma^{\infty}Y_{0}
}
\]
commutes in the category of Omega weak spectra.  According to
Lemma~\ref{lemspectrawkspectra}.(ii), we can choose a map of spectra
$\spapp{X}\to \spapp{Y}$ whose underlying map of Omega weak spectra is
the composite
\[
\spapp{X} \iso \csp{X}\overto{f} \csp{Y} \iso \spapp{Y}.
\]
Likewise, since $\cofapp{Y}\to \zspapp{Y}$ is an acyclic fibration, we
can choose a map of $B$-modules $\cofapp{X}\to \cofapp{Y}$ making the
diagram
\[
\xymatrix@R-=3ex{
\cofapp{X}\ar[r]\ar[d]&\zspapp{X}\ar[d]\\
\cofapp{Y}\ar[r]&\zspapp{Y}
}
\]
commute in $\Mod{B}$.  Then the following diagram commutes in $\owspcn$:
\[ \def\objectstyle{\scriptstyle}\def\labelstyle{\scriptscriptstyle}
\xymatrix@R-=3ex@C-=2.5ex{
\csp{X}\ar[d]_{f}
&\uw\spapp{X}\ar[l]\ar[r]\ar[d]
&\uw\stab\spapp{X}\ar[d]
&\uw\stab(\csp\bNCom \Sigma^{\infty}\cofapp{X})\ar[l]\ar[r]\ar[d]
&\uw\stab(\csp{Z}\Sigma^{\infty}\zspapp{X})\ar[d]
&\uw\csp{Z}\Sigma^{\infty}\zspapp{X}\ar[l]\ar[r]\ar[d]
&\csp{Z}\Sigma^{\infty}X_{0}\ar[d]^{f_{0}}\\
\csp{Y}
&\uw\spapp{Y}\ar[l]\ar[r]
&\uw\stab\spapp{Y}
&\uw\stab(\csp\bNCom \Sigma^{\infty}\cofapp{Y})\ar[l]\ar[r]
&\uw\stab(\csp{Z}\Sigma^{\infty}\zspapp{Y})
&\uw\csp{Z}\Sigma^{\infty}\zspapp{Y}\ar[l]\ar[r]
&\csp{Z}\Sigma^{\infty}Y_{0}
}
\]
The required naturality now follows.  (The same argument applied to
the identity map on $\csp{X}$ shows that $\bigiso{X}$ is in fact
independent of the choices.)  This completes the proof of
Theorem~\ref{intmainthm}. 

\section{Proof of Theorem~\ref{inthomthm}}\label{sechom}

In this section we prove Theorem~\ref{inthomthm}.  By
Proposition~\ref{propreduced}, it suffices to consider the case when
$B$ is a cofibrant commutative $R$-algebra, and by
Proposition~\ref{propki} it suffices to prove the analogous
theorem for reduced homology theories on the category $\NUCA{B}$.  
The work of Section~\ref{secspectra} reduces this to proving the
following theorem:

\begin{thm}\label{thmhom}
The category of reduced homology theories on $\NUCA{B}$ is equivalent
to the category of reduced homology theories on $\spc{\NUCA{B}}$.
\end{thm}

In this equivalence, the functor in one direction sends the reduced
homology theory $k_{*}$ in $\spc{\NUCA{B}}$ to the theory
$\restrict{k}{*}$ on $\NUCA{B}$ defined by 
\[
\restrict{k}{*}(N) = k_{*}(\dcsps{N}).
\]
where $\dcsps{N}$ denotes the left derived functor of the suspension
spectrum functor $\csps{}$.
On the other hand, for a reduced homology theory $h_{*}$ on
$\NUCA{B}$, we define a functor $\extend{h}{*}$ from $\spc{\NUCA{B}}$ to
graded abelian groups by 
\[
\extend{h}{q}(\csp{X}) = \Colim h_{q+n} X_{n}.
\]

\begin{lem}\label{lemc}
For any $\csp{X}$ in $\spc{\NUCA{B}}$, the map $\extend{h}{*}(\csp{X})\to
\extend{h}{*}(\stab\csp{X})$ is an isomorphism.
\end{lem}

\begin{proof}
The map $\extend{h}{q}\csp{X}\to
\extend{h}{q}\stab\csp{X}$ is the map
\begin{multline*}
\Colim_{n} h_{q+n} X_{n} \to
\Colim_{n} h_{q+n} (\Tel_{j\geq n} \Omega^{j-n} X_{j})\\
\iso
\Colim_{n} \Colim_{j\geq n} h_{q+n}(\Omega^{j-n}X_{j}).
\end{multline*}
To see that this is an isomorphism, it suffices to show that the map 
\[
\Colim_{n} h_{q+n} X_{n} \to
\Colim_{n} h_{q+n} (\Omega^{i}X_{n+i})
\]
is an isomorphism for all $i$.  The maps
\[
h_{q+n} (\Omega^{i}X_{n+i}) \iso
h_{q+n+i}(\dsusp{B}{}^{i} \Omega^{i}X_{n+i}) \to
h_{q+n+i} X_{n+i}
\]
induce a map
\[
\Colim_{n} h_{q+n} (\Omega^{i}X_{n+i})\to
\Colim_{n} h_{q+n} X_{n} 
\]
inverse to the map above.
\end{proof}

It follows that $\extend{h}{*}$ sends stable equivalences to isomorphisms
and therefore induces a functor from the stable category $\h\spc{\NUCA{B}}$ to the category of
graded abelian groups.  The suspension isomorphism for $h_{*}$ induces
a suspension isomorphism for $\extend{h}{*}$, and the Direct Sum Axiom for
$h_{*}$ implies the Direct Sum Axiom for $\extend{h}{*}$.  This then
defines a functor from the category of homology theories on $\NUCA{B}$
to the category of homology theories on $\spc{\NUCA{B}}$

In order to prove Theorem~\ref{thmhom}, it suffices to prove that
these functors $\extend{(-)}{}$ and $\restrict{(-)}{}$ are inverse
equivalences.  It is clear that when we start with a homology theory
$h_{*}$ on $\NUCA{B}$, we have a natural isomorphism between $h_{*}$
and the composite functor $\extres{h}{*}$.  To produce the natural
isomorphism between the other composite and the identity, we use the
following construction:

\begin{defn}
For a spectrum $\csp{X}$, let $\cutoff{n}{X}$ be the spectrum with
$j$-th object $X_{j}$ for $j\leq n$ and $\susp{B}^{j-n}X_{n}$ for $j>n$.
\end{defn}

We have compatible natural maps $\cutoff{n}{X}\to \cutoff{n+1}{X}$, with
$\csp{X}$ the colimit.  More usefully for the work below, the map
$\Tel \cutoff{n}{X}\to \csp{X}$ is an objectwise weak equivalence (a
weak equivalence on each object); this implies the following
less precise result: 

\begin{prop}
The natural map $\Tel \cutoff{n}{X}\to \csp{X}$ is a stable equivalence.
\end{prop}

If we write $\spfree{n}X_{n}$ for the spectrum that has $j$-th object
the initial object $*$ for $j<n$ and $\susp{B}^{j-n}X_{n}$ for $j\geq
n$, then we have natural maps 
\[
\susp{B}^{n}\cutoff{n}{X} \from
\susp{B}^{n}\spfree{n}X_{n} \to 
\csps{X_{n}}
\]
that are isomorphisms on $j$-th objects for $j\geq n$ and are therefore 
stable equivalences.  Applying $k_{*}$ and the suspension isomorphism,
we obtain the following proposition:

\begin{prop}
For $\csp{X}$ cofibrant, $k_{q}\cutoff{n}{X}\iso k_{q+n}\csps{X_{n}}$.
\end{prop}

The previous two propositions give us a natural (in both $\csp{X}$ and
$k_{*}$) isomorphism 
\[
\resext{k}{q}(\csp{X})=
\Colim k_{q+n}(\csps{X_{n}}) \to k_{q}(\csp{X}).
\]
Since this natural isomorphism commutes with the suspension isomorphisms, it
is a natural isomorphism of homology theories.  This completes the
proof of Theorem~\ref{thmhom}.

\section{Proof of Theorem~\ref{intdeloopthm}}\label{secdeloop}

In this section we prove Theorem~\ref{intdeloopthm} that interprets
the cotangent complex of the suspension spectra of \einf spaces in
terms of the associated spectra.  As indicated in the introduction, up
to equivalence, the definition of cotangent complex should not depend
on the \einf operad involved.  In order to take advantage of the work
in previous sections, we work with the linear isometries operad $\oL$,
and consider the category of $\oL$-spaces, the \einf spaces for the
\einf operad $\oL$.  If $X$ is an $\oL$-space, then
$\Sigma^{\infty}X_{+}$ is an ``$\oL$-spectrum'', an \einf ring
spectrum for the \einf operad $\oL$.  The category of $\oL$-spectra is
closely related to the category of commutative $S$-algebras: The
functor $S\sma_{\oL}(-)$ studied in EKMM \cite[I\S8]{ekmm} converts
$\oL$-spectra to weakly equivalent commutative $S$-algebras.  We study
the functor that takes an $\oL$-space $X$ to the 
commutative $S$-algebra $\newss{X}=S\sma_{\oL}\Sigma^{\infty}X_{+}$
and we prove the following theorem.

\begin{thm}\label{thmoldeloopthm}
For an $\oL$-space $X$, the $S$-module $\Scotan{\newss{X}}\iso
\derindec{S}\derI \newss{X}$ is
naturally isomorphic in the stable category to the spectrum associated
to $X$. 
\end{thm}

In the theorem, we are regarding $\newss{X}$ as augmented over
$S=\newss{*}$ via the map induced by the trivial map $X\to *$, and
$\Scotan{}$ is as in Section~\ref{secnuca} (with $A=B=S$).  
Since the cotangent complex of $\newss{X}$ is weakly equivalent to
the extended $\newss{X}$-module $(\newss{X})\sma \Scotan{\newss{X}}$,
Theorem~\ref{intdeloopthm} is an immediate consequence.

We analyze $\Scotan{\newss{X}}$ using the results of
Section~\ref{secspectra}.  Since the work of that section is phrased
in terms of model structures and Quillen adjunctions, it is convenient
to make a number of model category observations for $\oL$-spaces.
Since the category of $\oL$-spaces is the category of algebras for a
continuous monad, Quillen's small object argument and standard
techniques prove the following proposition.

\begin{prop}
The category of $\oL$-spaces is a topological closed model category
with weak equivalences and fibrations the weak equivalences and
(Serre) fibrations of the underlying spaces.
\end{prop}

May, Quinn, and Ray \cite[IV.1.8]{mayeinf} observe that
the functor $\Sigma^{\infty}(-)_{+}$ from $\oL$-spaces to
$\oL$-spectra is left adjoint to the zero-th space functor
$\Omega^{\infty}$.  Since the functor $\Omega^{\infty}$ preserves weak
equivalences and fibrations, this is in fact a Quillen adjunction.
Since the functor $S\sma_{\oL}(-)$ from
$\oL$-spectra to commutative $S$-algebras is a Quillen left adjoint,
we obtain the following result.

\begin{prop}
The functor $\newss{}$ from the category of
$\oL$-spaces to the category of commutative $S$-algebras is a
Quillen left adjoint.
\end{prop}

The previous propositions in particular give us a notion of cofibrant
$\oL$-space and prove that the suspension spectrum functor $\newss{}$
takes cofibrant $\oL$-spaces to cofibrant commutative $S$-algebras.
Since in $\oL$-spaces, the one-point $\oL$-space is both the initial and final
object, the category of $\oL$-spaces is enriched over the
category of based spaces.  When we regard $\newss{}$ as a functor into
the category of commutative $S$-algebras lying over $S$, the functor
$\newss{}$ is enriched over based spaces.  As a formal consequence we
obtain the following proposition.

\begin{prop}
The functor $\newss{}$ from the category of $\oL$-spaces to the
category of commutative $S$-algebras over $S$ preserves the tensor
with based spaces.  In particular, it converts the suspension functor
$\lsusp$ in the category of $\oL$-spaces to the suspension functor
$\susp{S}$ in the category of commutative $S$-algebras over $S$.
\end{prop}

We use the notation $\lsusp$ for the suspension in $\oL$-spaces
because of the following theorem, which is well-known to experts.  The
theorem is closely related to the uniqueness theorem of May and
Thomason \cite{maythomason} and the relationship between the homotopy
category of $\oL$-spaces and the homotopy category of connective
spectra. For the statement, note that since the category of
$\oL$-spaces may be defined as the category of algebras for a
continuous monad in
based spaces, the loop space 
$\Omega X$ is the underlying based space of 
cotensor of an $\oL$-space $X$ with the based space
$S^{1}$.

\begin{thm}\label{mtthm}
The derived functor of $\lsusp$ is a (one-fold) delooping functor:
If $X$ is a cofibrant $\oL$-space, then the unit of the suspension,
loops adjunction, 
\[
X\to \Omega \lsusp X
\]
is group completion.
\end{thm}

Since no proof of this theorem has appeared in the literature, we
outline a proof at the end of the section.  

The canonical maps $\Sigma \lsusp^{n}X\to \lsusp^{n+1}X$
make $\{\lsusp^{n}X\}$ a spectrum in the category of spaces, or in the
terminology of Lewis--May \cite{lms}, an ``indexed prespectrum''.
When $X$ is cofibrant, this has the further property that the 
adjoint structure map $\lsusp^{n}X\to \Omega \lsusp^{n+1}X$ is a weak
equivalence for $n>0$ and is group completion for $n=0$.
Also when $X$ is cofibrant, the structure maps are cofibrations, and so 
\[
Z=\Colim S^{-n} \sma \lsusp^{n} X,
\]
is a Lewis--May spectrum whose zeroth space is a group completion of $X$.  
Although this construction does not constitute an infinite loop space
machine on $\oL$-spaces, we have the following proposition;
see Remark~\ref{machine} below for further discussion.

\begin{prop}
If $X$ is cofibrant, $Z$ is a model for the spectrum associated to $X$.
\end{prop}

The Lewis--May spectrum $Z$ is naturally isomorphic in the stable
category to the $S$-module 
\[
Z_{S}=\Colim S^{-n}_{S} \sma \lsusp^{n} X.
\]
This is the composite of the ``free $\bL$-spectrum'' functor of EKMM
applied to $Z$, 
$\bL Z=\oL(1)\thp Z$,  and the functor $S\sma_{\oL}(-)$ from
$\bL$-spectra to $S$-modules.  For the proof of
Theorem~\ref{intdeloopthm}, it is useful to reframe this in the
context of the spectra in the category of $S$-modules, i.e., the
$\Mod{S}$-spectra of Section~\ref{secspectra}.  We have 
a $\Mod{S}$-spectrum $\csp{Z}$ defined by
\[
Z_{n}= S_{S} \sma \lsusp^{n}X,
\]
and we have a canonical natural map of $\Mod{S}$-spectra from
$\csp{Z}$ to the suspension $\Mod{S}$-spectrum of $Z_{S}$ induced by
the inclusion of 
$S_{S}\sma \lsusp^{n}X\iso \Sigma^{n}(S^{-n}_{S}\sma \lsusp^{n} X)$
in the colimit system defining $Z_{S}$.  An easy colimit argument
shows that this map is a stable equivalence.

The $\Mod{S}$-spectrum $\csp{Z}$ is one $\Mod{S}$-spectrum associated to
$\{\lsusp^{n}X\}$, but we also have a different one that
takes into account the action of the topological monoid $\oL(1)$ on
the based spaces $\lsusp^{n}X$.  The (Lewis--May) suspension spectrum functor
$\Sigma^{\infty}$ takes based $\oL(1)$-spaces to $\bL$-spectra; we
write $\newssp$ for the composite functor
$S\sma_{\oL}\Sigma^{\infty}(-)$ that lands in $S$-modules.  The
purpose for introducing this construction is that
we have a canonical natural isomorphism of $S$-modules 
\[
\newssp (X_{+}) \iso \newss{X}.
\]
We have a $\Mod{S}$-spectrum $\csp{B}$ defined by
\[
B_{n} = \newssp \lsusp^{n}X.
\]
The identity isomorphism of $\Sigma^{\infty}\lsusp^{n}X$ in the
category of Lewis--May spectra induces a map of $\bL$-spectra 
\[
(\bL S)\sma (\lsusp^{n}X) \iso \bL(\Sigma^{\infty}\lsusp^{n}X)\to 
\Sigma^{\infty}\lsusp^{n}X
\]
and a map of $S$-modules 
\[
S_{S}\sma \lsusp^{n}X = S\sma_{\oL}\bL S\sma \lsusp^{n}X \to 
S\sma_{\oL}\Sigma^{\infty}\lsusp^{n}X =
\newssp \lsusp^{n}X, 
\]
that induces a map of $\Mod{S}$-spectra $\csp{Z}\to \csp{B}$.
Since each map displayed above is a weak equivalence, the map
$\csp{Z}\to \csp{B}$  
stable equivalence.   We now have everything
needed for the proof of Theorem~\ref{thmoldeloopthm}.

\begin{proof}[Proof of Theorem~\ref{thmoldeloopthm}] 
It suffices to consider the case when $X$ is a cofibrant $\oL$-space
and show that $\derindec{S}\derI\newss{X}$ is naturally isomorphic to
$Z_{S}$ in the homotopy category of $S$-modules.  Applying
Propositions~\ref{propszmod} and~\ref{propindecforg}, and combining
with the work above, it suffices to show that the underlying
$\Mod{S}$-spectrum of the derived suspension spectrum
$\dcsps\derI\newss{X}$ is naturally isomorphic to $\csp{B}$ in the
homotopy category of $\spc{\Mod{S}}$.  Since $I$ is part of a Quillen
equivalence, its derived functor preserves suspension and
$\dcspsS\derI\newss{X}$ is naturally isomorphic to $\csp{I}^{\dR}
\scspsS \newss{X}$ in the homotopy category of $\spc{\Mod{S}}$. Since
the underlying $\Mod{S}$-spectrum of $\csp{I}^{\dR} \scspsS \newss{X}$
is the homotopy fiber of the augmentation, it is naturally weak
equivalent to the $\Mod{S}$-spectrum cofiber of the unit map:
\[
\scspsS \newss{X}\cup_{\scspsS S}C\scspsS S \overto{\sim}
\csp{I}^{\dR}\scspsS \newss{X}.
\]
Since the unit map is a cofibration, the cofiber is equivalent to the
quotient, $(\scspsS \newss{X})/ \scspsS S$.
Using the natural isomorphisms of $S$-modules
\[
(\susp{S}^{n}\newss{X})/S \iso (\newss{\lsusp^{n}X})/S
\iso \newssp{(\lsusp^{n}X_{+}/S^{0})}\iso \newssp{\lsusp^{n}X},
\]
we obtain our chain of natural isomorphisms in $\h\spc{\Mod{S}}$
between $\dcsps\derI\newss{X}$ and $\csp{B}$ as the zigzag
\[
\csp{I}^{\dR}\scspsS \newss{X} \overfrom{\sim}
\scspsS \newss{X}\cup_{\scspsS S}C\scspsS S \overto{\sim}
(\susp{S}^{n}\newss{X})/S \iso \newssp{\lsusp^{n}X} = \csp{B}.
\]
\end{proof}

We now go on to the proof of Theorem~\ref{mtthm}.  We write $\lfree$
for the monad on based spaces associated to the operad $\oL$.  For a
based space $T$, $\lfree T$ is the quotient of the disjoint union of
$\oL(n)\times_{\Sigma_{n}}T^{n}$ (cartesian power of $T$) by an
equivalence relation in terms of the basepoint and the operad
degeneracies (operadic multiplications with $\oL(0)=*$), described in
detail in \cite[2.4]{maygils}.

We study the suspension in $\oL$-spaces in terms of geometric
realization.  For an $\oL$-space $X$ and a based simplicial set
$T\subdot$, let $X\btensor {T\subdot}$ be the simplicial $\oL$-space
which in degree $n$ is the tensor of $X$ with the based set $T_{n}$;
this is the coproduct of copies of $X$ indexed on the non-basepoint
simplexes of $T_{n}$.  Writing $S^{1}\subdot$ for the simplicial
model of the based circle with one vertex and one non-degenerate
$1$-simplex, the following lemma implies in particular that $\lsusp X$
is the geometric realization of the simplicial $\oL$-space $X\btensor
S^{1}\subdot$.

\begin{lem}\label{lembar}
The geometric realization of a simplicial $\oL$-space is naturally an
$\oL$-space.  For any $\oL$-space $X$ and any based simplicial set
$T\subdot$, the map $X\btensor |T\subdot|\to |X\btensor T\subdot|$
induced by the universal property of the tensor is an isomorphism of
$\oL$-spaces. 
\end{lem}

\begin{proof}
The first statement is \cite[12.2]{maygils}:
Since cartesian products and colimits of spaces commute with geometric
realization, for any simplicial based space $Y\subdot$, we have a
natural isomorphism $\lfree |Y\subdot| \iso |\lfree Y\subdot|$. It is
straight-forward to check that the composite of this isomorphism and
the geometric realization of the $\oL$-space structure map $|\lfree
Y\subdot|\to |Y\subdot|$ provides an $\oL$-space structure map for
$|Y\subdot|$.  

For the statement about tensors, consider first the
free $\oL$-space $\lfree X$. The universal property of 
the free functor and the coproduct induce an isomorphism of simplicial 
$\oL$-spaces 
\[
\lfree(X\sma T\subdot)\iso (\lfree X)\btensor {T\subdot},
\]
and applying the isomorphism of the previous paragraph, we get an
isomorphism  
\[
(\lfree X)\btensor |T\subdot| \iso
\lfree(X\sma |T\subdot|)\iso \lfree|X\sma T\subdot|\iso 
|\lfree(X\sma T\subdot)| \iso |(\lfree X)\btensor{T\subdot}|,
\]
which is the statement for $\lfree X$. In the general case,
the tensor $X\btensor |T\subdot|$ is constructed as the reflexive
coequalizer 
\[
\xymatrix{
\lfree(\lfree X\sma |T\subdot|)\ar@<.75ex>[r]\ar@<-.75ex>[r]
&\lfree(X\sma |T\subdot|)\ar[r]&X\btensor |T\subdot|.
}
\]
Commuting geometric realization with $\lfree$ and the coequalizer, we see
that the coequalizer displayed above is the geometric realization of
the reflexive coequalizer 
\[
\xymatrix{
\lfree(\lfree X \sma T\subdot)
\ar@<.75ex>[r]\ar@<-.75ex>[r]
&\lfree(X\sma T\subdot)\ar[r]&X\btensor {T\subdot}
}
\]
describing the tensor of the $\oL$-space $X$ with the based set
$T\subdot$.
\end{proof}

The tensor of the $\oL$-space $X$ with a finite based set
defines a functor from finite based sets to based spaces that takes
the trivial based set $*$ to the trivial based set $*$.  This
constructs a $\Gamma$-space associated to $X$, and the previous lemma
identifies the suspension $\lsusp X$ as the classifying space of this
$\Gamma$-space.  Segal \cite{segalgamma} proved that when a
$\Gamma$-space is ``special'', the loop
space of the classifying space is a group completion.  In this case,
special means that the
map from $X\amalg\dots \amalg  X\to X\times \dots \times X$ is a weak
equivalence. Theorem~\ref{mtthm} is therefore an immediate consequence
of the following lemma.

\begin{lem}\label{lemcoprodprod}
If $X$ and $Y$ are cofibrant $\oL$-spaces, the map from the coproduct
$X\amalg Y$ to the cartesian product $X\times Y$ is a weak equivalence.
\end{lem}

We prove this using a shortcut.  Recall from \cite[9.6]{maygils}, the
two-sided monadic bar construction $B(\lfree,\lfree,X)$ which is the
geometric realization of the simplicial $\oL$-space
\[
B_{n}= \lfree \underbrace{\lfree \cdots \lfree}_{n} X.
\]
The iterated structure map $\lfree \cdots \lfree X\to X$ induces a map
of $\oL$-spaces $B(\lfree,\lfree,X)\to X$.  When $X$ is cofibrant, we can
choose a map of $\oL$-spaces $X\to B(\lfree,\lfree,X)$ and a homotopy
from the composite $X\to X$ to the identity through maps of $\oL$-spaces.
Choosing such a map for cofibrant $Y$ as well, we obtain a diagram
\[
\xymatrix{
X\amalg Y\ar[r]\ar[d]
&B(\lfree,\lfree,X)\amalg B(\lfree,\lfree,Y)\ar[r]\ar[d]
&X\amalg Y\ar[d]\\
X\times Y\ar[r]
&B(\lfree,\lfree,X)\times B(\lfree,\lfree,Y)\ar[r]
&X\times Y
}
\]
where the horizontal composites are homotopic to the identity (through
maps of $\oL$-spaces).  The lemma therefore reduces to showing that
the map 
\[
B(\lfree,\lfree,X)\amalg B(\lfree,\lfree,Y)\to
B(\lfree,\lfree,X)\times  B(\lfree,\lfree,Y)
\]
is a weak equivalence.  Both the coproduct of $\oL$-spaces and the
cartesian product commute with geometric realization.  Since this bar
construction is the geometric realization of a proper simplicial
space, we are reduced to showing that the map
\[
B_{n}(\lfree,\lfree,X)\amalg B_{n}(\lfree,\lfree,Y)\to
B_{n}(\lfree,\lfree,X)\times  B_{n}(\lfree,\lfree,Y)
\]
is a weak equivalence for all $n$.  The following lemma
therefore completes the proof of Lemma~\ref{lemcoprodprod}.

\begin{lem}\label{lastlem}
If $T$ and $U$ are nondegenerately based, then the map
$\lfree(T\vee U)\to \lfree T\times \lfree U$ is a homotopy equivalence.
\end{lem}

As always, ``nondegenerately based'' means that the inclusion of the
base point is an unbased $h$-cofibration.  The proof of
Lemma~\ref{lastlem}  
involves studying the double filtration on $\lfree(T\vee U)$ and 
$\lfree T\times \lfree U$ of homogeneous degree in $T$ and $U$: Let
$F^{m,n}\lfree(T\vee U)\subset \lfree(T\vee U)$ be the image of
$\oL(m+n)\times T^{m}\times U^{n}$, and let $F^{m,n}(\lfree T\times
\lfree U)\subset \lfree T\times \lfree U$ be $F^{m}\lfree T\times
F^{n}\lfree U$ where $F^{m}\lfree T\subset \lfree T$ is the image of
$\oL(m)\times T^{m}$, and similarly for $F^{n}\lfree U$. 
The map in Lemma~\ref{lastlem} preserves this double filtration.
Let $W^{m,n}$ denote the subspace of $T^{m}\times U^{n}$ consisting of
those points where at most $m+n-1$ coordinates are not the basepoint
(equivalently when $m,n>0$, the subset where at least one coordinate is the
basepoint).  Since we have assumed that $T$ and $U$ are
nondegenerately based, the inclusion of $W^{m,n}$ in $T^{m}\times
U^{m}$ is a $\Sigma_{m}\times \Sigma_{n}$-equivariant $h$-cofibration;
moreover, $F^{m,n}\lfree(T\vee U)$ is formed from
\[
F^{m-1,n}\lfree(T\vee U)\cup_{F^{m-1,n-1}\lfree(T\vee U)}
F^{m,n-1}\lfree(T\vee U)
\]
as the pushout over the map
\[
\oL(m+n)\times_{\Sigma_{m}\times \Sigma_{n}}W^{m,n}\to
\oL(m+n)\times_{\Sigma_{m}\times \Sigma_{n}}(T^{m}\times U^{n}).
\]
Likewise, $F^{m,n}(\lfree T\times \lfree U)$
is formed from
\[
F^{m-1,n}(\lfree T\times \lfree U)
\cup_{F^{m-1,n-1}(\lfree T\times \lfree U)}
F^{m,n-1}(\lfree T\times \lfree U)
\]
as the pushout over the map
\[
(\oL(m)\times \oL(n))\times_{\Sigma_{m}\times \Sigma_{n}}W^{m,n}\to
(\oL(m)\times \oL(n))\times_{\Sigma_{m}\times \Sigma_{n}}(T^{m}\times U^{n}).
\]
The maps
\begin{gather*}
\oL(m+n)\times_{\Sigma_{m}\times \Sigma_{n}}W^{m,n}\to
(\oL(m)\times \oL(n))\times_{\Sigma_{m}\times \Sigma_{n}}W^{m,n}\\
\noalign{\noindent and}
\oL(m+n)\times_{\Sigma_{m}\times \Sigma_{n}}(T^{m}\times U^{n})\to
(\oL(m)\times \oL(n))\times_{\Sigma_{m}\times \Sigma_{n}}(T^{m}\times U^{n})
\end{gather*}
are homotopy equivalences since the maps
$\oL(m+n)\to \oL(m)\times \oL(n)$ are equivariant homotopy
equivalences.  Since the map
\[
F^{0,0}\lfree(X\vee Y)\to F^{0,0}(\lfree T\times \lfree U)
\]
is the isomorphism $\oL(0)\to
\oL(0)\times \oL(0)$, an easy double induction shows that we have a
homotopy equivalence on $F^{m,n}$ for all $m,n$.  Passing to the
colimit, we see that the map $\lfree(T\vee U)\to\lfree T\times \lfree
U$ is a homotopy equivalence.  This completes the proof of
Lemma~\ref{lastlem}.

\begin{rem}\label{machine}
Let $X$ be a cofibrant $\oL$-space, and let $Z$ be as above.  This
remark explains in the terminology of May and Thomason
\cite{maythomason}, why $Z$ is equivalent to the output of an
``infinite loop space machine''.  Let $E$ be any infinite loop space
machine for $\oL$-spaces. Writing $\mathbf{n}$ for the finite based
set $\{0,\ldots,n\}$ (with $0$ as basepoint), the collection
$\{X\btensor \mathbf{n}\}$ forms a $\Gamma$-space in the category of
$\oL$-spaces, or an $\oF\oL$-space \cite[3.1]{maythomason}.  Applying
the machine $E$, a ``whiskering functor'', if necessary, and Segal's
machine, we obtain a bispectrum \cite[3.9ff]{maythomason}, that is
equivalent to $Z$ in one direction and $EX$ in the other.
\end{rem}

\section{Proof of Theorem~\ref{thmstable}}\label{secpfspectra}

This section is devoted to the proof of Theorem~\ref{thmstable}, the
topological closed model structure on the categories of spectra in
$\Com{B}{B}$, $\NUCA{B}$, and $\Mod{B}$.   For convenience we repeat
the definition of the cofibrations, fibrations, and weak equivalences;
we say that a map of spectra $\csp{X}\to \csp{Y}$ is:
\begin{enumerate}
\item A cofibration if $X_{0}\to Y_{0}$
is a cofibration and each $\susp{}Y_{n}\cup_{\susp{} X_{n}}X_{n+1}\to
Y_{n+1}$ is a
cofibration,
\item A fibration if each $X_{n}\to
Y_{n}$ is a fibration and each $X_{n}\to
Y_{n}\times_{\loops{}Y_{n+1}}\loops{}X_{n+1}$ 
is a weak equivalence, and
\item A stable equivalence if it induces an isomorphism of homotopy groups
$\pi_{*}\csp{X}\to\pi_{*}\csp{Y}$, where 
$\pi_{q}\csp{X}=\Colim_{n} \tpi_{q+n}X_{n}$
for $\tpi_{*}X=\pi_{*}X$ in $\NUCA{B}$ and $\Mod{B}$, and
$\tpi_{*}X=\Ker(\pi_{*}X\to\pi_{*}B)$ in $\Com{B}{B}$.
\end{enumerate}

It is clear that the categories of spectra in $\Com{B}{B}$,
$\NUCA{B}$, and $\Mod{B}$ have all small limits and
colimits. Likewise, it is clear from the definitions above that
cofibrations, fibration, and weak equivalences in these spectra are
closed under retracts and that weak equivalences have the
two-out-of-three property.  The proof of Theorem~\ref{thmstable}
therefore amounts to proving the factorization and lifting properties,
and proving the topological version of SM7.  The arguments are
identical for all of the categories, and we use $\CAT$ to denote any
of the categories $\Com{B}{B}$, $\NUCA{B}$, and $\Mod{B}$ in what
follows.

We begin with an alternative characterization of the acyclic fibrations.

\begin{lem}\label{lemaf}
A map $\csp{X}\to \csp{Y}$ is an acyclic fibration
in $\spc{\CAT}$ if and only if each map $X_{n}\to Y_{n}$ is an acyclic
fibration in $\CAT$.
\end{lem}

\begin{proof}
Since $\loops{}$ preserves fibrations and acyclic fibrations (even on
non-fibrant objects), when each $X_{n}\to Y_{n}$ is
an acyclic fibration in $\CAT$, the map $\csp{X}\to\csp{Y}$ is a
fibration and stable equivalence in $\spc{\CAT}$.
Conversely, assume $\csp{X}\to \csp{Y}$ is a fibration and stable
equivalence; then each map $X_{n}\to Y_{n}$ is a fibration, and so it
suffices to show that each map $X_{n}\to 
Y_{n}$ is a weak equivalence.  Let $\csp{W}$ be the fiber (for
$\CAT=\Com{B}{B}$, this means $W_{n}=B\times_{Y_{n}}X_{n}$); then
$\csp{W}$ is an Omega spectrum.  Since 
the sequential colimit (of abelian groups) is an exact functor, the
levelwise long 
exact sequences of homotopy groups for $W_{n}\to X_{n}\to Y_{n}$
induce a long exact sequence of 
homotopy groups for $\csp{W}\to \csp{X}\to \csp{Y}$.  It follows that
$\pi_{*}W=0$ and therefore that $\tpi_{*}W_{n}=0$ for all $n$, since $\csp{W}$ is
an Omega spectrum.  We see from the long exact sequence of homotopy groups for
$W_{n}\to X_{n}\to Y_{n}$ that $X_{n}\to Y_{n}$ is a weak equivalence. 
\end{proof}

The previous lemma allows us to prove the lifting property for
cofibrations and acyclic fibrations:  If the solid rectangle
\[
\xymatrix@R-1em{
\csp{A}\ar[r]\ar@{>->}[d]&\csp{X}\ar@{->>}[d]^{\sim}\\
\csp{B}\ar[r]\ar@{-->}[ur]&\csp{Y}
}
\]
is a commutative diagram with the left-hand map a cofibration and the
right-hand map an acyclic fibration of $\CAT$-spectra, we can construct
the dashed arrow making the diagram commute as follows.  We have that
$A_{0}\to B_{0}$ is a cofibration and $X_{0}\to Y_{0}$ is an acyclic
fibration, and so we use the lifting property in $\CAT$ to construct
the required map $B_{0}\to X_{0}$.  Inductively, having constructed
$B_{n}\to X_{n}$, since we have that
$\susp{}B_{n}\amalg_{\susp{}A_{n}}A_{n+1}\to B_{n+1}$ is a cofibration
and $X_{n+1}\to Y_{n+1}$ is an acyclic fibration, we can use the
lifting property in $\CAT$ to construct a map $B_{n+1}\to X_{n+1}$,
compatible with the map $\susp{}B_{n}\to\susp{}X_{n}\to X_{n+1}$ and
making the required diagram commute.  This proves the following proposition.

\begin{prop}\label{propcl}
In $\spc{\CAT}$, cofibrations have the left lifting property with
respect to acyclic fibrations.
\end{prop}

We could construct the factorization for cofibrations and acyclic
fibrations analogously, but in order to analyze the topological
version of SM7, it is useful to construct them instead using Quillen's
small object argument.  For this, we recall the set $I_{\CAT}$ of
``generating cofibrations'' and $J_{\CAT}$ of ``generating acyclic
cofibrations'' in $\CAT$, the definition of which is implicit in the
construction of the model structure on $\CAT$ in EKMM
\cite[VII\S4]{ekmm}.  For $\CAT=\Mod{B}$, $I$ is simply the set of
``cells''
\[
I_{\Mod{B}}=\{ S^{m}_{B}\to CS^{m}_{B}\mid m\in \bZ \},
\]
where
$S^{m}_{B}$ denotes the cofibrant sphere $B$-module
\cite[III\S2]{ekmm} and $J$ is the set of cylinders of spheres,
\[
J_{\Mod{B}}=\{ S^{m}_{B}\to S^{m}_{B}\sma I_{+}\mid m\in \bZ \}.
\]
For $\CAT=\NUCA{B}$, the sets $I$ and $J$ are 
\begin{gather*}
I_{\NUCA{B}}=\bNCom I_{\Mod{B}} = \{\bNCom i\mid i\in I_{\Mod{B}}\}\\
J_{\NUCA{B}}=\bNCom J_{\Mod{B}} = \{\bNCom j\mid j\in J_{\Mod{B}}\}.
\end{gather*}
For $\CAT=\Com{B}{B}$, the sets $I$ and $J$ are only slightly more
complicated: $I$ is the set of diagrams of commutative $B$-algebras
\[
\xymatrix@R-1em{
\bCom S^{m}_{B}\ar[rr]\ar[dr]&&\bCom CS^{m}_{B}\ar[dl]\\
&B
}
\]
where $\bCom$ denotes the free commutative $B$-algebra, and the map
$S^{m}_{B}\to CS^{m}_{B}$ is always the usual inclusion.  The set $J$ has an
entirely analogous description using the inclusion $S^{m}_{B}\to
S^{m}_{B}\sma I_{+}$.  The fundamental property of the sets $I_{\CAT}$
and $J_{\CAT}$ is the following: 

\begin{prop}\label{propij}
A map in $\CAT$ is an acyclic fibration if and only if it has the right lifting
property with respect to the maps in $I_{\CAT}$; it is a fibration if
and only if it has the right lifting 
property with respect to the maps in $J_{\CAT}$.
\end{prop}

In order to describe sets of generating cofibrations and generating
acyclic cofibrations for $\spc{\CAT}$, we need one more piece of
notation.  For an object $X$ in $\CAT$, we let $\spfree{n}X$ denote
the $\CAT$-spectrum with $(\spfree{n}{X})_{j}$ the initial object for
$j<n$ and $\susp{}^{j-n}X$ for $j\geq n$; maps of $\CAT$-spectra from
$\spfree{n}X$ into a $\CAT$-spectrum $\csp{Y}$ are in one-to-one
correspondence with maps in $\CAT$ from $X$ to $Y_{n}$.  We set $\ISP$
to be the set of maps
\[
\ISP=\{\spfree{n}f\mid f
\in I_{\CAT}, n=0,1,2,\dots  \}.
\]
We have the following analogue of the first part of the previous
proposition. 

\begin{lem}\label{lemisp}
A map in $\spc{\CAT}$ is an acyclic fibration if and only if it has
the right lifting property with respect to $\ISP$.
\end{lem}

The description of the set $\JSP$ is slightly more complicated.
Certainly $\JSP$ should 
contain the maps $\spfree{n}f$ for $f\in J_{\CAT}$, but these maps
only generate the cofibrations that are levelwise weak equivalences.
Another general sort of stable equivalence occurs in the following
way: If $f\colon S\to T$ is a map in $I_{\CAT}$, then
for any map $S\to X_{n}$, we could attach the cell $\spfree{n}f$ at
the $n$-th level 
or attach the cell $\spfree{n+1}\susp{}f$ at the $(n+1)$-st level, and
the map
\[
\csp{X}\amalg_{\spfree{n+1}(\susp{}S)}\spfree{n+1}(\susp{}T) \to
\csp{X}\amalg_{(\spfree{n}S)}\spfree{n}T
\]
is a stable equivalence.  This map
is not a cofibration, but we can make it a cofibration using a
cylinder:  If we denote by $(-)\otimes I$ the tensor in $\CAT$ of
an object with the (unbased) unit interval, the map
\begin{multline*}
\csp{X}\amalg_{\spfree{n+1}(\susp{}S)}
\spfree{n+1}(\susp{} S \otimes I)\amalg_{\spfree{n+1}(\susp{}S)}
\spfree{n+1}(\susp{}T) \\\to
(\csp{X}\amalg_{\spfree{n}S}\spfree{n}T)\amalg_{\spfree{n+1}(\susp{}T)}
\spfree{n+1}(\susp{}T\otimes I)
\end{multline*}
gives a version of the previous map that is a cofibration.
With this as motivation, for $f\colon S\to T$ in $I_{\CAT}$, let
\begin{align*}
\srs{n}{f}&=\spfree{n}S\amalg_{\spfree{n+1}(\susp{}S)}
\spfree{n+1}(\susp{} S \otimes I)\amalg_{\spfree{n+1}(\susp{}S)}
\spfree{n}(\susp{}T),\\
\tar{n}{f}&=\spfree{n}T\amalg_{\spfree{n+1}(\susp{}T)}
\spfree{n+1}(\susp{} T \otimes I),
\end{align*}
and let $\jmap{n}{f}\colon \srs{n}{f}\to \tar{n}{f}$ be the map
induced by $f$.  Then on $j$-th objects, $\jmap{n}{f}$ is the identity
map (on $*$) for $j<n$, is the map $S\to T$ (from $I_{\CAT}$) for
$j=n$, and is the map
\[
\susp{}^{j-n}S\otimes I\amalg_{\susp{}^{j-n}S}\susp{}^{j-n}T
\to \susp{}^{j-n}T\otimes I,
\]
for $j>n$.  This last map is easily seen to be the inclusion of a
deformation retraction.  (In fact, it is isomorphic to a map in $J_{\CAT}$.)

In particular, we have that $\jmap{n}{f}$ is a cofibration and for
$j>n$ is a weak equivalence on $j$-th objects.  It follows that
$\jmap{n}{f}$ is an acyclic cofibration.
We set $\JSP$ to be the set of maps
\[
\JSP = \{\spfree{n}g\mid g
\in J_{\CAT}, n=0,1,2,\dots  \}
\cup
\{ \jmap{n}{f}
\mid f\in I_{\CAT}, n=0,1,2,\dots  
\}
\]

In studying these maps, it is convenient to use the following
notation: Let
\[
M_{n}\csp{X} =
X_{n}\times_{\loops{} X_{n+1}}
(\loops{} X_{n+1})^{I},
\]
where $(-)^{I}$ denotes the cotensor in $\CAT$ with the unbased
interval.  The two endpoint of the interval induce two maps $(\loops{}
X_{n+1})^{I}\to \loops{} X_{n+1}$.  The construction of $M_{n}\csp{X}$
uses one of these maps; the other gives us a map $M_{n}\csp{X}\to
\loops{} X_{n+1}$.  Unwinding the definition of of the maps
$\jmap{n}{f}$ and the universal property of tensors and cotensors
leads to the following proposition.

\begin{prop}\label{propjmap}
Let $f\colon S\to T$ be a map in $I_{\CAT}$ and let $\csp{X}$ be a
$\CAT$-spectrum.  Then maps in $\spc{\CAT}$ from
$\tar{n}{f}$ to $\csp{X}$ are in one-to-one correspondence with maps
in $\CAT$ from $T$ to $M_{n}\csp{X}$.  Maps in $\spc{\CAT}$ from
$\srs{n}{f}$ to $\csp{X}$ are in one-to-one correspondence with
commutative diagrams in $\CAT$,
\[
\xymatrix{
S\ar[d]_{f}\ar[r]&M_{n}\csp{X}\ar[d]\\
T\ar[r]&\loops{} X_{n+1}.
}
\]
\end{prop}

We need one more observation about the map $M_{n}\csp{X}\to \loops{}
X_{n+1}$ before moving on to the analogue for $\JSP$ of
Lemma~\ref{lemisp}.

\begin{lem}\label{lemhalfjsp}
Let $\csp{X}\to \csp{Y}$ be a map in $\spc{\CAT}$ and assume that
the maps $X_{n}\to Y_{n}$ are fibrations in $\CAT$ for all $n$.  Then the map
\[
M_{n}\csp{X}\to M_{n}\csp{Y}\times_{\loops{} Y_{n+1}}\loops{} X_{n+1}
\]
is a fibration in $\CAT$.
\end{lem}

\begin{proof}
Abbreviate $M_{n}\csp{X}$ to $M$ and $M_{n}\csp{Y}\times_{\loops{}
Y_{n+1}}\loops{} X_{n+1}$ to $N$; we have a commutative cube where
the double-headed arrows are known to be fibrations and the dotted
arrow is the map we want to show is a fibration.
\[
\xymatrix@C-2pc{
M\ar[rr]\ar@{->>}[dd]\ar@{..>}[dr]
&&(\loops{} X_{n+1})^{I}\ar@{->>}[dr]\ar@{->>}'[d][dd]\\
&N\ar@{->>}[dd]\ar[rr]
&&(\loops{} Y_{n+1})^{I}\ar@{->>}[dd]\\
X_{n}\times \loops{} X_{n+1}\ar'[r][rr]\ar@{->>}[dr]
&&\loops{} X_{n+1}\times \loops{} X_{n+1}\ar@{->>}[dr]\\
&Y_{n}\times \loops{} X_{n+1}\ar[rr]
&&\loops{} Y_{n+1}\times \loops{} Y_{n+1}
}
\]
The front and back (rectangular) faces are pullbacks, and the map
\[
(\loops{} X_{n+1})^{I}
\to 
(\loops{} Y_{n+1})^{I}
\times_{(\loops{} Y_{n+1}\times \loops{} Y_{n+1})}
\loops{} X_{n+1}\times \loops{} X_{n+1}
\]
is a fibration.  Since fibrations in $\CAT$ are characterized by the
right lifting property with respect to $I_{\CAT}$, it follows that
$M\to N$ is a fibration.
\end{proof}

\begin{lem}\label{lemjsp}
A map in $\spc{\CAT}$ is a fibration if and only if it has
the right lifting property with respect to $\JSP$.
\end{lem}

\begin{proof}
Given $h\colon \csp{X}\to \csp{Y}$, it follows from Proposition~\ref{propij}
that the maps $X_{n}\to Y_{n}$ are fibrations for all $n$ if and only
if $h$ has the right lifting property with respect to the set
$\{\spfree{n}g\mid g
\in J_{\CAT}, n=0,1,2,\dots  \}$.
We can therefore restrict to the case when $X_{n}\to Y_{n}$ is a
fibration for all $n$ and prove that the map
$X_{n}\to Y_{n}\times_{\loops{}Y_{n+1}}\loops{}X_{n+1}$ 
is a weak equivalence if and only if $h$ has the right lifting
property with respect to
$\{ \jmap{n}{f}
\mid f\in I_{\CAT}, n=0,1,2,\dots  
\}$.  
Proposition~\ref{propjmap} 
implies that $h$ having the right lifting property with respect to
the maps $\jmap{n}{f}$ for $f\in I_{\CAT}$ (for fixed $n$) is equivalent
to 
\[
M_{n}\csp{X}\to M_{n}\csp{Y}\times_{\loops{} Y_{n+1}}\loops{} X_{n+1}
\]
having the right lifting property with respect to $I_{\CAT}$, which is
equivalent to it being an acyclic fibration (by
Proposition~\ref{propij}).  By the previous lemma, the displayed map
is a fibration, so being an acyclic fibration is equivalent to being a weak
equivalence. 
\end{proof}

Recall that for a set of maps $A$ (e.g., $A=\ISP$ or $A=\JSP$), a
relative $A$-complex is a map $\csp{X}\to \Colim \csp{X}_{n}$ where
$\csp{X}_{0}=\csp{X}$ and each $\csp{X}_{n+1}$ is formed from
$\csp{X}_{n}$ as the pushout over a coproduct of maps in $A$.
Lemmas~\ref{lemisp} and~\ref{lemjsp} therefore give us the right lifting
property of acyclic fibrations and fibrations with respect to
relative $\ISP$-complexes and relative $\JSP$-complexes
(respectively).  We have already observed that the maps in $\ISP$ and
$\JSP$ are cofibrations and it follows that relative $\ISP$-complexes
and relative $\JSP$-complexes are cofibrations.  The remainder of the
proof of following lemma requires a
compactness argument that we give at the end of the section.

\begin{lem}\label{acxcof}
A relative $\ISP$-complex is a cofibration.  A relative $\JSP$-complex
is an acyclic cofibration.
\end{lem}

The following lemma constructs factorizations:

\begin{lem}\label{acxfact}
Let $A=\ISP$ or $\JSP$.  Any map $\csp{X}\to \csp{Y}$ can be factored
as a relative $A$-complex
$\csp{X}\to \csp{Z}$ and a map $\csp{Z}\to \csp{Y}$ that has the right
lifting property with respect to the maps in $A$.
\end{lem}

\begin{proof}
Proposition~\ref{propjmap}  (for $A=\JSP$) and the characterization of maps
out of $\spfree{n}$ (for $A=\ISP$) show that for $\csp{A}$ the
domain or codomain of a map in $A$, the set of maps of $\CAT$-spectra
out of $\csp{A}$ commutes with sequential colimits,
\[
\Colim \spc{\CAT}(\csp{A},\csp{X}_{n}) 
\iso \spc{\CAT}(\csp{A},\Colim \csp{X}_{n}),
\]
when the maps $\csp{X}_{n}\to \csp{X}_{n+1}$ are cofibrations.
We can now apply Quillen's small object argument to construct the
required factorizations.
\end{proof}

The usual retract argument (factoring using the previous lemma and
applying the lifting property of Proposition~\ref{propcl}) then proves
the following lemma, the converse of Lemma~\ref{acxcof}.

\begin{lem}\label{lemcharc}
If a map in $\spc{\CAT}$ is a cofibration, then it is a retract of a
relative $\ISP$-complex; if it is an acyclic cofibration, then it is a
retract of a relative $\JSP$-complex.
\end{lem}

We have now assembled everything we need for the proof of the theorem.

\begin{proof}[Proof of Theorem~\ref{thmstable}]
As observed above, the proof that classes of maps defined in the
statement form a closed model structure is completed by proving the
required factorization and lifting properties.
Using the characterization of the acyclic
cofibrations from \ref{lemcharc}, the lifting
properties follow from Proposition~\ref{propcl} and
Lemma~\ref{lemjsp}.  The factorization properties follow from
Lemma~\ref{acxcof} and Lemma~\ref{acxfact}.

It remains to prove the topological version of SM7: We need to show
that when $i\colon
\csp{A}\to \csp{B}$ is a cofibration and $p\colon \csp{X}\to \csp{Y}$
is a fibration, the map of spaces
\[
\spc{\CAT}(\csp{B},\csp{X})\to
\spc{\CAT}(\csp{B},\csp{Y})
\times_{\spc{\CAT}(\csp{A},\csp{Y})}
\spc{\CAT}(\csp{A},\csp{X})
\]
is a (Serre) fibration and is a weak equivalence if either $i$ or $p$ is a
stable equivalence.  

To show that the map is a fibration, it suffices to consider the case
when $i$ is a relative $\ISP$-complex by Lemma~\ref{lemcharc}, and
for this, it suffices to consider the case when $i$ is a map in
$\ISP$.  Then $i$ is a map $\spfree{n}{f}\colon \spfree{n}{S}\to \spfree{n}{T}$
for some $n$ and some $f$ in $I_{\CAT}$, and we can identify the map in
question with the map of spaces
\[
\CAT(T,X_{n})\to
\CAT(T,Y_{n})
\times_{\CAT(S,Y_{n})}
\CAT(S,X_{n}).
\]
This is a fibration by the topological version of SM7 for $\CAT$.  An
entirely similar argument proves that this map is an acyclic fibration
when $p$ is an acyclic fibration.

Finally, we need to show that the map is an acyclic fibration when $i$
is an acyclic cofibration.  As in the previous paragraph, this reduces
to the case when $i$ in in $\JSP$.  When $i$ is $\spfree{n}g$ for some
$g$ in $J_{\CAT}$, the argument reduces to $\CAT$ just as in the
previous paragraph.  Now consider the other case, when $i=\jmap{n}{f}$ for some $f$
in $I_{\CAT}$. The argument is (as always) to go back over the proof
of the lifting property of Lemma~\ref{lemjsp} taking into account the
topology of the mapping spaces.  
Proposition~\ref{propjmap} was stated in terms of
a bijection of sets, but the argument refines to give an isomorphism
of spaces; this allows us to identify the map in question with the map
of spaces 
\[
\CAT(T,M)\to \CAT(T,N)\times_{\CAT(S,N)}\CAT(S,N),
\]
where we have used the notation in the proof of
Lemma~\ref{lemhalfjsp}.  Lemma~\ref{lemjsp} shows that when
$\csp{X}\to \csp{Y}$ is a fibration, $M\to N$ is an acyclic fibration.
It now follows from the
topological version of SM7 in $\CAT$ that 
the map displayed above is an acyclic fibration.
\end{proof}

Finally, we complete the proof of Lemma~\ref{acxcof} by proving that a
relative $\JSP$-complex is a stable equivalence.  For this, 
it suffices to see that a
pushout 
\[
\csp{Y}=\csp{X}\amalg_{(\coprod \srs{n_{\alpha}}{f_{\alpha}})}
(\coprod \tar{n_{\alpha}}{f_{\alpha}})
\]
over a coproduct of maps $\jmap{n_{\alpha}}{f_{\alpha}}$ is a stable
equivalence.  This is 
clear for a finite coproduct (since then for $n$ large, the map on
$n$-th objects is a deformation retract).  The proof for the general
case follows from the finite case, provided we can identify the
homotopy groups $\pi_{*}\csp{Y}$ as the filtered colimit of the
homotopy groups $\pi_{*}\csp{Y}_{A}$ where $A$ ranges over the finite
subsets of the index sets.  For this it is sufficient to identify the
homotopy groups of the $n$-th object of $\csp{Y}$ as the filtered
colimit of the homotopy groups of the $n$-th objects.  As observed
above, for each $\alpha$, $\jmap{n_{\alpha}}{f_{\alpha}}$ is on $n$-th
objects either an isomorphism (if $n<n_{\alpha}$), the map
$f_{\alpha}$ from $I_{\CAT}$ (when $n=n_{\alpha}$), or the inclusion
of a deformation retraction (if $n>n_{\alpha}$).  The argument is
therefore completed by the following lemma.

\begin{lem}\label{lemcompact}
Let $\{f_{\alpha }\colon S_{\alpha}\to T_{\alpha }\}$ be a set of maps
in $I_{\CAT}$, let $X$ be an object in $\CAT$ and let 
\[
Y=X\amalg_{(\coprod S_{\alpha})}(\coprod T_{\alpha}).
\]
Any element of $\tpi_{q}Y$ is
represented in $\tpi_{q}Y_{\alpha_{1},\dots,\alpha_{m}}$, 
\[
Y_{\alpha_{1},\dots,\alpha_{m}}=
X\amalg_{(S_{\alpha_{1}}\amalg \dots \amalg S_{\alpha_{m}})}
(T_{\alpha_{1}}\amalg \dots \amalg T_{\alpha_{m}}),
\]
for some finite subset of indexes $\alpha_{1},\dots,\alpha_{m}$.  If
some element of $\tpi_{q}Y_{\alpha_{1},\dots,\alpha_{m}}$ is zero in
$\tpi_{q}Y$, then it is zero in
$\tpi_{q}Y_{\alpha_{1},\dots,\alpha_{p}}$ for some subset of indexes
$\alpha_{1},\dots,\alpha_{p}$.
\end{lem}

\begin{proof}
In order to fix notation, we treat the case $\CAT=\Com{B}{B}$, but the
arguments for the other cases are similar. (We continue to write
$\amalg$ instead of $\sma$ for the coproduct, however, to avoid
confusion for the other cases.) Since
$\tpi_{q}$ is a subset of $\pi_{q}$, it suffices to prove the
analogous lemma for $\pi_{q}$, and for this, it suffices to work in
$\ComB$, the category of commutative $B$-algebras.  In this context,
each map $f_{\alpha}$ is just a map $\bCom S^{m_{\alpha}}_{B}\to \bCom
CS^{m_{\alpha}}_{B}$ for some integer $m_{\alpha}$.  For a finite
subset $A=\{\alpha_{1},\dots,\alpha_{n}\}$, we set
$Y_{A}=Y_{\alpha_{1},\dots,\alpha_{n}}$ (or $Y_{A}=X$ when $A$ is
empty), and we set 
\[
M_{A}=\bigvee_{\alpha \notin A}S^{m_{\alpha }}_{B},
\]
so that $Y\iso Y_{A}\amalg_{\bCom M_{A}}\bCom(CM_{A})$.
We construct various filtrations (of $B$-modules) on $Y$ using the bar
construction, which we denote as ``$\beta\subdot$''.  Let 
$\beta\subdot(Y_{A},\bCom M_{A},\bCom *)$ be the simplicial object in
$\ComB$, with 
\[
\beta_{n}(Y_{A},\bCom M_{A},\bCom*) = Y_{A} \amalg 
\underbrace{\bCom M_{A}\amalg \dots \amalg \bCom M_{A}}_{n} 
{}\amalg
\bCom *,
\]
with faces induced by the maps $M_{A}\to Y$, $M_{A}\to *$, and the
codiagonal maps, and with degeneracy maps induced by inserting
an extra summand of $\bCom M_{A}$.  We write $\beta[A]$ for the
geometric realization and we use the notation
\[
Y_{A}=\beta[A]_{0}\to \beta[A]_{1}\to \cdots \to \beta[A]_{n}\to\cdots,
\qquad \beta[A]=\Colim \beta[A]_{n},
\]
for the filtration arising from the geometric realization.
Since the degeneracy maps are inclusions
of wedge summands of $B$-modules, this filtration is 
a filtration by $h$-cofibrations of
$B$-modules.  The geometric realization of a simplicial object in $\ComB$
is an object in $\ComB$, and we have isomorphisms
\begin{multline*}
\beta[A]=|\beta\subdot(Y_{A},\bCom M_{A},\bCom *)| 
\iso 
Y_{A}\amalg_{\bCom M_{A}} 
|\beta\subdot(\bCom M_{A},\bCom M_{A},\bCom M_{A})|
\amalg_{\bCom M_{A}}\bCom *\\
\iso Y_{A} \amalg_{\bCom M_{A}} 
\bCom(M_{A}\sma I)
\amalg_{\bCom M_{A}}\bCom *
\iso Y_{A}\amalg_{\bCom M_{A}}\bCom (CM_{A})\iso Y
\end{multline*}
(see \cite[VII.3.2]{ekmm}).  For $A\subset A'$, the map $\beta[A]\to
\beta[A']$ covering the identity map of $Y$ is not induced by a
simplicial map but does preserve the filtrations above.

Now given an element $x$ of $\pi_{q}Y$, we show that $x$ is the image
of an element of $\pi_{q}Y_{A}$
for some finite subset of indexes $A$.  Let $A_{0}$ be the empty set.
Using the isomorphisms 
\[
\pi_{q}Y\iso \pi_{q}(\beta[A_{0}])\iso \Colim \pi_{q}(\beta[A_{0}]_{n}),
\]
we can represent $x$ as the image of an element $x_{0}$ of
$\pi_{q}(\beta[A_{0}]_{n})$ for some $n$.  Now suppose by induction,
we have constructed a finite set $A_{k}$ and found an element $x_{k}$
of $\pi_{q}(\beta[A_{k}]_{n-k})$ whose image in $\pi_{q}Y$ is $x$.
The quotient $\beta[A_{k}]_{n-k}/\beta[A_{k}]_{n-k-1}$ is the quotient
of $\Sigma^{n-k}\beta_{n-k}(Y_{A_{k}},\bCom M_{A_{k}},\bCom *)$ by the
degeneracies.  Since the degeneracies are inclusions of wedge
summands, this quotient is still a wedge sum of $B$-modules, each
summand of which involves only finitely many of the indexes $\alpha$.
The image of $x_{k}$ factors through $\pi_{q}$ of a finite wedge sum
of these summands; let $A_{k+1}$ be the union of $A_{k}$ and the
finite set of indexes involved in these summands.  Since the map
$\beta[A_{k}]\to \beta[A_{k+1}]$ is compatible with the filtration,
and by construction, the image of $x_{k}$ in
$\pi_{q}(\beta[A_{k+1}]_{n-k}/\beta[A_{k+1}]_{n-k-1})$ is zero, we can
find an element $x_{k+1}$ in $\pi_{q}(\beta[A_{k+1}]_{n-k-1})$ whose
image in $\pi_{q}(\beta[A_{k+1}]_{n-k})$ is the image of $x_{k}$.  It
follows that the image of $x_{k+1}$ in $\pi_{q}Y$ is $x$.  Continuing in this
way, we get a finite subset of indexes $A_{n}$ and an element $x_{n}$
of $\pi_{q}(\beta[A_{n}]_{0})=\pi_{q}Y_{A_{n}}$ whose image in
$\pi_{q}Y$ is $x$.

The argument for a relation is similar, using a relative
class in $\pi_{q+1}$ and starting with
$A_{0}=\{\alpha_{1},\ldots,\alpha_{m}\}$. 
\end{proof}

\section{Cohomology Theories for Operadic Algebras}\label{secoperadic}

We have stated and proved the results in this paper in terms of the
special case of particular interest, the category of algebras over the
operad $\oCom$.  Most of these results hold quite generally for
the categories of algebras over other operads.  The purpose of this
section is to give precise statements of these general results and to
indicate how to adapt the arguments in the earlier sections to the
more general case.  We prove the following theorem.

\begin{thm}\label{thmgenop}
Let $\oGen$ be an operad of (unbased) spaces with each $\oGen(n)$ of
the homotopy type of a $\Sigma_{n}$-CW complex.  Let $B$ be a
cofibrant $\oGen$-algebra in EKMM $R$-modules, and let $UB$ be its
universal enveloping algebra.  Let $\aGen$ be the category of
$\oGen$-algebras of EKMM $R$-modules lying over $B$.
\begin{enumerate}
\renewcommand{\theenumi}{\arabic{enumi}}\renewcommand{\labelenumi}{\theenumi.}
\item Topological Quillen Cohomology with coefficients in a
$UB$-module induces an equivalence from
the homotopy category of left $UB$-modules to the category of cohomology
theories on $\aGen$.
\item The category of homology theories on $\aGen$ is equivalent to
the category of homology theories on the category of left $UB$-modules.
\item The stable category of $\oGen$-algebras over and under $B$ is
equivalent to the category 
of left $UB$-modules.
\item The equivalence in the previous statement takes the suspension
spectrum of an $\oGen$-algebra $A$ over $B$ to the $UB$-module of
infinitesimal $UB$-deforma\-tions.
\end{enumerate}
\end{thm}

We review the general definition of the universal enveloping algebra
$UB$ and of Topological Quillen Cohomology below; the proof of the
theorem essentially amounts to formulating the definitions in a
framework parallel to the case for $\oGen=\oCom$.  In the case
$\oGen=\oCom$, $UB$ is just $B$, the category of left $UB$-modules is
the category of $B$-modules, and Topological Quillen Cohomology is
Topological Andr\'e--Quillen Cohomology, as in the theorems in the
introduction.   

In general, unlike the case of the operad $\oCom$, a weak equivalence
of $\oGen$-algebras does not necessarily induce a weak equivalence of
enveloping algebras, and this is why we need to assume that $B$ is
cofibrant from the outset in the theorem above.  The theorem combined
with Proposition~\ref{propreduced} implies that for general $B$, the
category of cohomology theories on $\aGen$ is equivalent to the
homotopy category of left $UB'$-modules and the category of homology
theories on $\aGen$ is equivalent to the category of homology theories
on the category of left $UB'$-modules, where $B'\to B$ is a cofibrant
approximation.

For $\oGen=\oAss$, the operad for associative algebras, $UA$ is
$A\sma_{R} A^{\op}$, and so the category of left $UA$-modules is the
category of $A$-bimodules.  When $A$ is cofibrant, one typically
writes $A^{e}$ for $A\sma_{R} A^{\op}$. More generally, $A^{e}$ denotes
$A'\sma_{R} A^{\prime\op}$, for some fixed choice of cofibrant
approximation $A'\to A$. 
Lazarev \cite{lazarev} identifies Topological Quillen Cohomology in
terms of Topological Hochschild Cohomology, and identifies the module
of infinitesimal deformations of an associative algebra $A$ as the
homotopy fiber of the multiplication map $A^{e}\to A$. Part 1 of the
previous theorem then has the following corollary.

\begin{cor}
Let $B$ be an associative $R$-algebra.  Every cohomology
theory on the category of associative $R$-algebras lying over $B$ is
of the form
\[
h^{*}(X,A) = \pi_{-*}Fib(THH_{R}(X,M) \to THH_{R}(A,M))
\]
for some left $B^{e}$-module $M$.
\end{cor}

We now return to the general case of Theorem~\ref{thmgenop}.  We begin
by describing the universal enveloping algebra $UB$.  At the same time, 
we describe the ``universal enveloping operad'' $\oUn$.

\begin{defn}
The universal enveloping operad $\oUn$ is the operad in $R$-modules
that has $n$-th object $\oUn(n)$ defined by the coequalizer 
\[
\xymatrix{
\bigvee_{k} \oGen(n+k)\sma_{\Sigma_{k}}(\bGen B)^{(k)}
\ar@<-.75ex>[r]\ar@<.75ex>[r]
&\bigvee_{k} \oGen(n+k)\sma_{\Sigma_{k}}B^{(k)}
\ar[r]
&\oUn(n)
}
\]
where $\bGen$ denotes the free $\oGen$-algebra functor, one map is the
$\oGen$-action map of $B$ and the other is induced by the operadic
multiplication of $\oGen$.  The operadic multiplication on $\oUn$ is
induced by the operadic multiplication of $\oGen$.  The universal
enveloping algebra $UB$ is the $R$-algebra $\oUn(1)$.
\end{defn}

The fundamental property of the universal enveloping operad is given
by the following proposition, which is an easy consequence of the definitions.

\begin{prop}\label{propuniversal}
The category of $\oGen$-algebras lying under $B$ is equivalent to the
category of $\oUn$-algebras.
\end{prop}

In particular, we have that $\oUn(0)=B$.  Let $\oNUn$ be the operad with 
\[
\oNUn(n)=
\begin{cases}
*&n=0\\
\oUn(n)&n>0
\end{cases}
\]
The category of $\oNUn$-algebras plays the role for
$\oGen$-algebras under and over $B$
that the category of nucas plays for commutative algebras.

\begin{defn}
Let $\aGenB$ denote the category of $\oUn$-algebras lying over $B$,
and let $\aGenN$ denote the category of $\oNUn$-algebras.  Let
$K\colon \aGenN\to \aGenB$ denote the functor that takes a
$\oNUn$-algebra $N$ to $B\vee N$.  Let $I\colon \aGenB\to\aGenN$
denote the functor that takes $A$ to the (point-set) fiber of the
augmentation $A\to B$.
\end{defn}

The functors $K$ and $I$ are adjoint, and the following proposition
holds just as in the commutative algebra case. 

\begin{thm}
The categories $\aGenB$ and $\aGenN$ are topological closed model
categories with weak equivalences the weak equivalences of the
underlying $R$-modules.  The adjunction $(K,I)$ is a Quillen
equivalence. 
\end{thm}

\begin{proof}
The topological closed model structures and easy consequences of the
general theory in EKMM 
\cite[VII\S4]{ekmm}.  When we take $\bT$ to be the monad associated to
the operad $\oUn$ or $\oNUn$, or indeed any operad in the category of $R$-modules, the
proof of the ``Cofibration Hypothesis''of  
\cite[VII\S4]{ekmm} for $\bT$ follows just like the proof for
associative and commutative algebras in 
\cite[VII\S3]{ekmm}. The
key observation is that since colimits and smash
products of $R$-modules commute with geometric realization of
simplicial $R$-modules, the monad $\bT$ commutes with geometric
realization.  This gives the geometric realization of a simplicial
$\bT$-algebra a $\bT$-algebra structure and also proves that geometric
realization commutes with colimits of $\bT$-algebras.  
As a consequence, we obtain the analogue of \cite[VII.3.7]{ekmm}: For
any maps of $\bT$-algebras $A\to A'$ and $A\to 
A''$,  we can identify the geometric
realization of the ``bar construction''
\[
\beta_{n} (A',A,A'')
=(A,\bT M,\bT*) = A' \amalg 
\underbrace{A\amalg \dots \amalg A}_{n} 
{}\amalg A'',
\]
as the double pushout in $\bT$-algebras 
\[
A'\amalg_{A}(A\otimes I)\amalg_{A}A'',
\]
where $A\otimes I$ is the tensor with the (unbased) interval (see for
example the argument for Lemma~\ref{lembar}).  In the
special case when $A''=\bT *$ and $A=\bT M$ for some $R$-module $M$,
the degeneracy maps are the inclusion of wedge summands, and so the
filtration on the geometric realization is a filtration by
$h$-cofibrations. In particular, in this case, the inclusion of the
lowest filtration level $A'=A'\amalg \bT *$ in the geometric
realization, $A'\amalg_{\bT M}\bT(CM)$, is an $h$-cofibration; this is
the Cofibration Hypothesis. 

Since $I$ preserves fibrations an acyclic
fibrations, $(K,I)$ is a Quillen adjunction.  Finally, given any
$\oNUn$-algebra $N$, and any fibrant $\oUn$-algebra $X$ over $B$, it
is clear from the effect of $K$ and $I$ on homotopy groups that a map
$N\to IX$ is a weak equivalence if and only if the adjoint map $KN\to
X$ is a weak equivalence, and so $(K,I)$ is a Quillen equivalence.
\end{proof}

The proof of the previous theorem used the free functor from
$R$-modules to $\oNUn$-algebras, but there is in addition a free functor
$\bNUn$ from left $UB$-modules to $\oNUn$-algebras, left adjoint to the
forgetful functor from $\oNUn$-algebras to left $UB$-modules, defined by
\[
\bNUn M = \bigvee_{n\geq 1} (\oUn(n) \sma_{UB^{(n)}}M^{(n)})/\Sigma_{n},
\]
for a left $UB$-module $M$.
We also have a zero
multiplication functor $Z$ that gives $M$ a $\oNUn$-action where 
\[
\oNUn(n)\sma_{UB^{(n)}} M^{(n)}\to M
\]
is the left $UB$-action map for $n=1$ and the trivial map for $n>1$.
This functor has a left adjoint $Q$ defined by the coequalizer 
\[
\xymatrix{
\bNUn X\ar@<-.75ex>[r]\ar@<.75ex>[r]&X\ar[r]&QX,
}
\]
where one map is the $\oNUn$-algebra action map on $X$, and the other
map is the unit on the $\oNUn(1)\sma_{UB}X$ summand and
the trivial map on the summands $(\oUn(n)
\sma_{UB^{(n)}}M^{(n)})/\Sigma_{n}$ for $n>1$.  Since the zero
multiplication functor $Z$ preserves fibrations and weak equivalences,
we obtain the following proposition.

\begin{prop}
The $(Q,Z)$ adjunction is a Quillen adjunction.
\end{prop}

As a consequence of the previous two propositions, for a cofibrant
$\oGen$-algebra $A$ over $B$ and any left $UB$-module $M$, we obtain
bijections of  
sets (in fact, isomorphisms of abelian groups) 
\begin{multline*}
\h\aGen(A,B\vee ZM)\iso
\h\aGenB(B\amalg A,B \vee ZM) \\\iso 
\h\aGenB(I^{\dR}(B\amalg A),ZM)\iso 
\h\Mod{UB}(Q^{\dL}I^{\dR}(B\amalg A),M),
\end{multline*}
where $\dL$ and $\dR$ denote left and right derived functors.
This leads to the following definition.

\begin{defn}
For a $\oNUn$-algebra $N$, the module of infinitesimal
$UB$-deforma\-tions is the left $UB$-module $Q^{\dL}N$.  For a
cofibrant $\oGen$-algebra $A$ over $B$, the $UB$-module of infinitesimal
$UB$-deforma\-tions of $A$ is the left $UB$-module
$Q^{\dL}I^{\dR}(B\amalg A)$, and for a cofibration of $\oGen$-algebras
$A\to X$ over $B$, the 
$UB$-module of infinitesimal
$UB$-deforma\-tions of $X$ relative to $A$ is the left $UB$-module
$Q^{\dL}I^{\dR}(B\amalg_{A}X)$.
For a left $UB$-module $M$, we define
the Topological Quillen Cohomology of $X$ relative to $A$ with
coefficients in $M$ by 
\[
D^{*}_{\oGen}(X,A;M)=\Ext^{*}_{UB}(Q^{\dL}I^{\dR}(B\amalg_{A}X),M).
\]
For general $A$ and a general map $A\to X$ of $\oGen$-algebras over
$B$, the Topological Quillen Cohomology is defined using a cofibration
$A'\to X'$ covering $A\to X$ for cofibrant
approximations $A'\to A$ and $X'\to X$.
\end{defn}

Since $(I,K)$ is a Quillen equivalence, $I^{\dR}$ preserves coproducts
and cofibration sequences.  Since $Q$ is a Quillen left adjoint, $Q^{\dL}$
also preserves coproducts and cofibration sequences.   From this, it
is easy to see that for any left $UB$-module $M$, Topological Quillen
Cohomology with coefficients in $M$ forms a cohomology theory, with
the connecting maps $\delta$ induced by $\Ext$.

We now have the underlying theory for $\oGen$-algebras parallel to the
theory for commutative algebras, and the argument for
Theorem~\ref{thmgenop} parallels the proof of
Theorems~\ref{intmainthm}--\ref{intcotanc}.  The Brown's
Representability argument in Section~\ref{secreduced} applies
generally to any category to which the arguments of EKMM
\cite[VII\S4]{ekmm} apply.  The only argument in
Sections~\ref{secreduced}--\ref{sechom} and~\ref{secpfspectra} that does
not immediately 
generalize in this framework is the appeal to \cite[III.5.1]{ekmm} in
Lemma~\ref{lemsymprod}. 
(For the proof of Lemma~\ref{lemcompact}, note that the coproduct of
$\oUn$-algebras $X\amalg \bUn M$ has as its underlying
module $\bigvee (\oUnX(n)\sma_{R} M^{(n)})/\Sigma_{n}$ by 
the analogue of Proposition~\ref{propuniversal} for $X$.)
The following lemma fills this gap.

\begin{lem}\label{lemgenone}
Let $M$ be a cofibrant left $UB$-module.  The natural map 
\[
E\Sigma_{m+}\sma_{\Sigma_{m}}(\oUn(m)\sma_{UB^{(m)}}M^{(m)})
\to (\oUn(m)\sma_{UB^{(m)}}M^{(m)})/\Sigma_{m}
\]
is a weak equivalence.
\end{lem}

The proof depends strongly on the hypothesis that the spaces of
$\oGen$ have equivariant CW homotopy types and that $B$ is a cofibrant
$\oGen$-algebra.  This latter hypothesis implies that $B$ is a retract
of a ``cell $\oGen$-algebra'' \cite[VII.4.11]{ekmm}.  If $B\to B'\to B$
is such a retraction (with $B\to B$ the identity), then we get a
retraction of operads $\oUn\to \oUn' \to \oUn$.
The analogue of Lemma~\ref{lemgenone} for $B'$,
then implies the lemma as stated for $B$.  Thus, it suffices to
consider the case when $B$ is a cell $\oGen$-algebra. Specifically,
this means that we can write $B$ as $\Colim B_{n}$ where 
$B_{0}=\oGen(0)_{+}\sma R$, and $B_{n+1}=B_{n}\amalg_{\bGen W_{n}}\bGen
(CW_{n})$, where $W$ is a wedge of sphere $R$-modules $S^{m}_{R}$.

The filtration of $B$ allows us to get a better hold
on the $\Sigma_{m}$-equivariant $R$-modules $\oUn(m)$.  For example,
$\oUn_{0}(m)=\oGen(m)_{+}\sma R$, and 
\[
\oUn_{1}(m)=\bigvee_{k\geq 0} \oGen(m+k)\sma_{\Sigma_{k}}(\Sigma W_{0})^{(k)}.
\]
More generally, we have the following lemma.

\begin{lem}\label{unfilt}
For each $n\geq 0$, $\oUn_{n+1}(m)$ has a filtration by
$\Sigma_{m}$-equivariant $h$-cofibrations
$\oUn_{n+1}=\Colim_{k} \oUn_{n+1}(m)_{k}$, with
$\oUn_{n+1}(m)_{0}=\oUn_{n}(m)$, and the 
filtration quotients
\[
\oUn_{n+1}(m)_{k}/\oUn_{n+1}(m)_{k-1}\iso 
\oUn_{n}(m+k)\sma_{\Sigma_{k}}(\Sigma W_{n})^{(k)}.
\]
\end{lem}

\begin{proof}
We set $\oUn_{n+1}(m)_{0}=\oUn_{n}(m)$ as required.
The idea is that $\oUn_{n+1}(m)_{k}$ is the image in $\oUn_{n+1}(m)$
of $\oUn_{n}(m+k)\sma_{\Sigma_{k}}(CW_{n})^{(k)}$.  Precisely,
consider the $\Sigma_{k}$-equivariant filtration
\[
W_{n}^{(k)}=F^{0}CW_{n}\to \cdots \to F^{k-1}CW_{n}\to
F^{k}CW_{n}=(CW_{n})^{(k)} 
\]
obtained as the smash power of the filtration on $CW_{n}$ that has
$W_{n}$ in level zero and $CW_{n}$ in level one.  For $k\geq 1$,
define $\oUn_{n+1}(m)_{k}$ as the pushout
\[
\xymatrix{
\oUn_{n}(m+k)\sma_{\Sigma_{k}}F^{k-1}CW_{n}\ar[r]\ar[d]
&\oUn_{n}(m+k)\sma_{\Sigma_{k}}(CW_{n})^{(k)}\ar[d]\\
\oUn_{n+1}(m)_{k-1}\ar[r]
&\oUn_{n+1}(m)_{k}.
}
\]
Then it is clear that the map
$\oUn_{n+1}(m)_{k-1}\to \oUn_{n+1}(m)_{k}$
is an $h$-cofibration, and the quotient is as indicated in the
statement.  The identification of $\oUn_{n+1}$ 
with $\Colim_{k}\oUn_{n+1}(m)_{k}$
follows from the universal properties and the fact that the map 
\[
\oUn_{n}(m+k)\sma_{1\times \Sigma_{k-1}}(W_{n}\sma (CW_{n})^{(k-1)})
\to
\oUn_{n}(m+k)\sma_{\Sigma_{k}}F^{k-1}CW_{n}
\]
is a categorical epimorphism.
\end{proof}

\begin{proof}[Proof of Lemma~\ref{lemgenone}]
By the usual retract argument, it suffices to prove the lemma when $M$
is a cell left $UB$-module, and the usual filtration argument then reduces
to the case when $M=UB\sma_{R} S^{q}_{R}$, that is, to proving that the
map
\[
E\Sigma_{m+}\sma_{\Sigma_{m}}(\oUn(m)\sma_{R}(S^{q}_{R})^{(m)})
\to (\oUn(m)\sma_{R}(S^{q}_{R})^{(m)})/\Sigma_{m}
\]
is a weak equivalence.  
This now follows from
previous lemma using the argument of
\cite[\S9]{mbthesis} (which generalizes the argument of
\cite[III.5.1]{ekmm}).  
\end{proof}

\section{Weak Equivalences and Excision in Model Categories}
\label{secpairs}

The axioms listed in the introduction for a homology or cohomology
theory on a closed model category $\CAT$ patently depend on the
cofibrations.  The reduced version of these axioms in
Section~\ref{secreduced} implicitly depends on the cofibrations in
terms of the definition of cofibration sequences.  The purpose of this
section is to prove the following theorem. 

\begin{thm}\label{thmpars}
Let $\CAT_{1}$ and $\CAT_{2}$ be closed model structures on the same category
with the same weak equivalences.  Let $h^{*}$ be a contravariant
functor from the category of pairs to the 
category of graded abelian groups and let $\delta^{n} \colon h^{n}(A)\to
h^{n+1}(X,A)$ be natural transformations of abelian groups.  Then
$(h^{*},\delta)$ is a cohomology theory on $\CAT_{1}$ if and only if
it is a cohomology theory on $\CAT_{2}$.  Likewise, for a covariant
functor and natural transformations, $(h_{*},\partial)$
is a homology theory on $\CAT_{1}$ if and only if it is a
homology theory on $\CAT_{2}$.
\end{thm}

As a technical point, with the Product Axiom and Direct Sum Axiom as
stated in the introduction, we need the standard assumption that
the closed model category has all small colimits. See
Remark~\ref{nocoprods} below for further discussion of the case when
this assumption does not hold.

The theorem above does not appear to be well known, and we offer it
here for its intrinsic interest; it has not been used in the previous
sections.  The basic idea is that although cofibrations determine the
notion of excision in a category, an appropriate notion of ``weak
excision'' defined in terms of ``homotopy cocartesian'' diagrams
depends only on the weak equivalences.  The following definition is
standard.

\begin{defn}
A commutative diagram 
\[
\xymatrix{
A\ar[r]\ar[d]&X\ar[d]\\
B\ar[r]&Y
}
\]
in a closed model category is homotopy cocartesian means that there
exists a commutative diagram
\[
\xymatrix{
X'\ar[d]^{\sim}
&A'\ar@{>->}[r]\ar@{>->}[l]\ar[d]^{\sim}
&B'\ar[d]^{\sim}\\
X&A\ar[l]\ar[r]&B
}
\]
with $A'$ cofibrant, 
the top horizontal arrows cofibrations, and vertical arrows weak
equivalences, such that the induced map $Y'=X'\cup_{A'}B'\to Y$ is a
weak equivalence. 
\end{defn}

More concisely, the diagram is homotopy cocartesian if the canonical
map in the homotopy category from the homotopy pushout to $Y$ is an
isomorphism.  One class of examples of homotopy cocartesian squares is
given by the squares where $A$ is cofibrant, $A\to X$ and $A\to B$ are
cofibrations, and $Y$ is the pushout $X\cup_{A}B$. Another class of
examples is given by the squares where $A\to B$ and $X\to Y$ are weak
equivalences.  With these examples in mind, the following proposition
is clear from the definition.

\begin{prop}
Let $\CAT$ be a closed model category and let $h$ be a contravariant
functor from the category of pairs to the category of abelian groups.
The following are equivalent:
\begin{itemize}\upshape
\item[(a)] $h$ satisfies the Homotopy Axiom and the Excision Axiom:
\begin{itemize}
\item[(i)](Homotopy) If $(X,A)\to (Y,B)$ is a weak equivalence of
pairs, then $h(Y,B)\to h(X,A)$ is an isomorphism. 
\item[(iii)](Excision) If $A$ is cofibrant, $A\to B$ and $A\to X$ are
cofibrations, and $Y$ is the pushout $X\cup_{A}B$, then the map of
pairs $(X,A)\to (Y,B)$ induces an isomorphism $h(Y,B)\to h(X,A)$.
\end{itemize}
\item[(b)] $h$ satisfies the following Weak Excision Axiom:  Whenever
\[
\xymatrix@R-1pc{
A\ar[r]\ar[d]&X\ar[d]\\
B\ar[r]&Y
}
\]
is homotopy cocartesian,
the map of pairs $(X,A)\to
(Y,B)$ induces an isomorphism $h(Y,B)\to h(X,A)$.
\end{itemize}
\end{prop}

The Weak Excision Axiom does not depend on the cofibrations but only
on the weak equivalences:

\begin{lem}
Let $\CAT_{1}$ and $\CAT_{2}$ be closed model structures on the same
category $\CAT$ with the same weak equivalences.  Then a commutative
diagram 
\[
\xymatrix@R-1pc{
A\ar[r]\ar[d]&X\ar[d]\\
B\ar[r]&Y
}
\]
is homotopy cocartesian in $\CAT_{1}$ if and only if it is homotopy
cocartesian in $\CAT_{2}$.
\end{lem}

\begin{proof}
Assume the diagram is homotopy cocartesian in $\CAT_{1}$.  Using the
factorization properties of $\CAT_{2}$ (starting with $A$), we can find
a commutative diagram
\[
\xymatrix@R-1pc{
X_{2}\ar[d]^{\sim}
&A_{2}\ar[d]^{\sim}\ar@{>->}[l]_{2}\ar@{>->}[r]^{2}
&B_{2}\ar[d]^{\sim}\\
X&A\ar[l]\ar[r]&B
}
\]
with $A_{2}$ cofibrant in $\CAT_{2}$, the top horizontal arrows
cofibrations in $\CAT_{2}$, and the vertical arrows weak
equivalences.  Likewise, using the
factorization properties of $\CAT_{1}$ and $\CAT_{2}$, we can extend
this to a commutative diagram
\[
\xymatrix@R-1pc{
X'\ar[d]^{\sim}
&A'\ar[d]^{\sim}\ar@{>->}[l]_{2}\ar@{>->}[r]^{2}
&B'\ar[d]^{\sim}\\
X_{1}\ar[d]^{\sim}
&A_{1}\ar[d]^{\sim}\ar@{>->}[l]_{1}\ar@{>->}[r]^{1}
&B_{1}\ar[d]^{\sim}\\
X_{2}\ar[d]^{\sim}
&A_{2}\ar[d]^{\sim}\ar@{>->}[l]_{2}\ar@{>->}[r]^{2}
&B_{2}\ar[d]^{\sim}\\
X&A\ar[l]\ar[r]&B
}
\]
where $A_{1}$ is cofibrant in $\CAT_{1}$, $A'$ is cofibrant in
$\CAT_{2}$, the horizontal arrows labeled with the number $i$ are
cofibrations in $\CAT_{i}$, and all the vertical arrows are weak
equivalences.  Taking $Y_{i}=X_{i}\cup_{A_{i}} B_{i}$ and
$Y'=X'\cup_{A'}B'$, the previous diagram induces maps
\[
Y' \to Y_{1}\to Y_{2}\to Y.
\]
Clearly the map $Y'\to Y_{2}$ is a weak equivalence (see, for example,
the characterization of homotopy pushouts by Dwyer and Spalinski
\cite[10.7]{modcat}). The hypothesis that the 
original diagram is homotopy cocartesian in $\CAT_{1}$ implies that the map
$Y_{1}\to Y$ is a weak equivalence.  It follows that the map $Y'\to
Y$ is a weak equivalence, and so the original diagram is homotopy
cocartesian in $\CAT_{2}$. 
\end{proof}

Since coproducts of cofibrant objects represent the coproduct in the
homotopy category, it is clear that the Product Axiom of the
introduction is equivalent to the following axiom.
\begin{itemize}
\item [$\text{(\ref{productaxiom})}_{w}$]
If $\{X_{\alpha}\}$ is a set of
objects and $X$ is the coproduct in the homotopy category, then the
natural map $h^{*}(X)\to \prod h^{*}(X_{\alpha})$ is an isomorphism. 
\end{itemize}
The analogous observation holds for the Direct Sum Axiom.
Since the homotopy category depends only on the weak equivalences in
the model structure and not the cofibrations, this completes the proof
of Theorem~\ref{thmpars}.

\begin{rem}\label{nocoprods}
It is sometimes useful to consider model categories that do not have
all small colimits but (as in the original definition in \cite{quil})
are only assumed to have finite colimits.  For these categories, it
appears unlikely that the version of the Product Axiom above
is equivalent to the one in the introduction, and it depends on the
application which axiom, if either, is the ``right'' one.  Typically,
the most useful version of the Product Axiom in this case is one where
we assume the isomorphism only for index sets of certain fixed
cardinalities; when coproducts of the given cardinalities always exist
in the point-set category, then the two versions of the axioms are
again equivalent.  The version of the Product Axiom for finite
cardinalities follows from the other axioms.
\end{rem}


\bibliographystyle{plain}

\end{document}